\documentclass{article}

\usepackage[utf8]{inputenc}
\usepackage[T1]{fontenc}
\usepackage{url}
\usepackage{booktabs}
\usepackage{nicefrac}
\usepackage{microtype}
\usepackage[dvipsnames,table]{xcolor}
\usepackage{amsmath,amsthm,amsfonts,amssymb,mathrsfs,float}
\usepackage{mathtools}
\usepackage{algorithm}
\usepackage{algorithmic}
\usepackage{comment}
\usepackage{bm}
\usepackage{bbm}
\usepackage{graphicx}
\usepackage{verbatim}
\usepackage{dsfont}
\usepackage{relsize}
\usepackage{caption}
\usepackage{subcaption}
\usepackage[english]{babel}
\usepackage[numbers,sort&compress]{natbib}

\usepackage{tabulary}
\usepackage{colortbl}
\usepackage{enumitem}
\usepackage{tablefootnote}
\usepackage{pifont}
\usepackage{geometry}
\usepackage{hyperref}
\hypersetup{
	colorlinks=true,
	linkcolor=red,
	anchorcolor=blue,
	citecolor=blue,
	filecolor=magenta,
	menucolor=red,
	urlcolor=cyan,
	hypertexnames=false,
	pdftitle={Incremental Aggregation on the Grassmannian for Asynchronous Eigenspace Computation},
	pdfauthor={Xiaolu Wang, Jiang Hu, Hoi-To Wai}
}
\usepackage[capitalize,noabbrev]{cleveref}

\newtheorem{assumption}{Assumption}
\newtheorem{lemma}{Lemma}
\newtheorem{proposition}{Proposition}
\newtheorem{theorem}{Theorem}
\newtheorem{corollary}{Corollary}
\newtheorem{remark}{Remark}

\crefname{assumption}{Assumption}{Assumptions}
\Crefname{assumption}{Assumption}{Assumptions}

\DeclareMathOperator*{\argmin}{arg\,min}
\DeclareMathOperator{\grad}{grad}

\DeclareMathOperator{\diag}{diag}
\DeclareMathOperator{\dist}{dist}
\DeclareMathOperator{\tr}{tr}
\DeclareMathOperator{\St}{St}
\DeclareMathOperator{\Gr}{Gr}

\DeclareMathOperator{\Log}{Log}
\DeclareMathOperator{\Polar}{Polar}

\newcommand{\comm}[2]{\mathcal{C}(#1,#2)}
\newcommand{\cmark}{\ding{51}}%
\newcommand{\xmark}{\ding{55}}%
\def\algname{GRASSIA}

\title{Incremental Aggregation on the Grassmannian for Asynchronous Eigenspace Computation}

\author{
	Xiaolu Wang\thanks{Software Engineering Institute, East China Normal University, Shanghai, China. Email: \href{mailto:xiaoluwang@sei.ecnu.edu.cn}{\texttt{xiaoluwang@sei.ecnu.edu.cn}}.}
	\and Jiang Hu\thanks{Yau Mathematical Sciences Center, Tsinghua University, Beijing, China. Email: \href{mailto:jianghu@tsinghua.edu.cn}{\texttt{jianghu@tsinghua.edu.cn}}.}
	\and Hoi-To Wai\thanks{Department of Systems Engineering and Engineering Management, The Chinese University of Hong Kong, Hong Kong SAR, China. Email: \href{mailto:htwai@se.cuhk.edu.hk}{\texttt{htwai@se.cuhk.edu.hk}}.}
}

\date{}

\begin{document}

\maketitle

\begin{abstract}
We study asynchronous optimization for finite-sum eigenspace computation in heterogeneous distributed systems. The theoretical foundations for asynchronous eigenspace computation remain scarce, with existing approaches offering limited coverage of dynamics directly on the Grassmannian under stale information. In this paper, we propose a Grassmannian incremental aggregation method that refreshes only arriving components and reuses cached gradients, retaining low per-update cost without global synchronization. The method employs an extrinsic polar update that preserves the intrinsic subspace geometry without requiring parallel transport of stale tangent vectors. Our analysis establishes a tight angle-dependent gradient-dominance characterization of the objective and a basin-invariance property for stale aggregated updates. These yield two-phase linear convergence, comprising an explicit broad-basin regime and a sharper local regime, with constants controlled by component spectral spreads. Experiments on serial and distributed PCA demonstrate improved sample efficiency and wall-clock convergence over representative baselines.
\end{abstract}


\section{Introduction}

We study asynchronous distributed computation of the leading-$k$ eigenspace of a symmetric matrix.
Eigenspace computation is a core primitive underpinning a broad range of applications, including principal component analysis (PCA) \citep{jolliffe2005principal}, singular value decomposition (SVD) \citep{golub1970singular}, low-rank approximation \citep{eckart1936approximation,halko2011finding}, and communication- or memory-efficient training schemes for large machine learning models \citep{vogels2019powersgd,zhao2024galore}.
Given an $n \times n$ symmetric matrix $\bm A\in\mathbb S^d$,
by the Ky Fan maximum principle, the leading-$k$ eigenspace problem can be written as the orthogonally constrained trace minimization problem \citep{golub2013matrix,saad2011numerical}
\begin{align}
	\min_{\bm{W} \in \St(d,k)} F(\bm{W}) \coloneqq -\tr(\bm{W}^\top \bm{A} \bm{W}),
	\label{eq:eigen}
\end{align}
where $\St(d,k)\coloneqq\{\bm W\in\mathbb R^{d\times k}:\bm W^\top\bm W=\bm I_k\}$ denotes the Stiefel manifold.
As the objective function is invariant under right orthogonal transformations $\bm W\mapsto\bm W\bm Q$, the decision variable represents a point on the Grassmannian rather than an ordered basis.
If the target eigengap $\delta\coloneqq\lambda_k(\bm A)-\lambda_{k+1}(\bm A)$ is positive, every global minimizer of \eqref{eq:eigen} spans the unique leading-$k$ invariant subspace.
In the distributed setting, we note that $\bm A$ is available through a \emph{finite-sum} representation,
\begin{align}
	\bm{A} \coloneqq \frac{1}{n}\sum_{i=1}^{n}\bm{A}_i,
	\qquad
	F(\bm{W}) = \frac{1}{n}\sum_{i=1}^{n} F_i(\bm{W}),
	\qquad
	F_i(\bm{W}) \coloneqq -\tr(\bm{W}^\top \bm{A}_i \bm{W}),
	\label{eq:finite-sum}
\end{align}
where the $i$th component ${\bm A}_i$ can be formed from a local data shard, a block, a time window, or a separately queried oracle, and is stored by the $i$th agent or worker in the system. 
As a common setting in parallel and distributed optimization \citep{bertsekas2015parallel,assran2020advances}, we aim to solve \eqref{eq:eigen} by updating the subspace iterate whenever a worker result arrives, without waiting for a synchronized round.

Algorithms for solving \eqref{eq:eigen} have a long history. Classical approaches include orthogonal iteration \citep{saad2011numerical}, Krylov methods such as Lanczos \citep{lanczos1950iteration}, Riemannian gradient methods \citep{absil2008optimization,alimisis2024geodesic}, and randomized subspace iteration \citep{halko2011finding}. These methods construct successive approximations through repeated applications of the global operator $\bm A$ and, under an eigengap and suitable initialization, typically converge
\emph{linearly} to the leading eigenspace.
Such algorithms can be viewed as \emph{synchronous} algorithms under \eqref{eq:finite-sum} as they require collecting all local contributions $\bm A_i\bm W$ before forming the next update. 
This leads to the \emph{straggler} issue in a distributed setting as some workers can be less efficient than others, thus resulting in significant synchronization overhead that limits the performance of solving \eqref{eq:eigen}. 
Distributed eigenspace methods reduce some of this burden through local computation, subspace averaging, or round-based aggregation \citep{huang2020communication,alimisis2021distributed,chen2021decentralized,li2021communication,grammenos2020federated,guo2024fedpower,ye2021deepca,gang2021distributed,andrade2023distributed}, but they still suffer from a round-based structure vulnerable to stragglers. 

An alternative solution to the straggler issue is to apply the incremental eigenspace methods which evaluate a single component $F_i$ at each iteration. Representative examples include Oja's method \citep{oja1982simplified,oja1985stochastic} and its variants \citep{mitliagkas2013memory,hardt2014noisy,huang2021streaming}, as well as Krasulina's method \citep{krasulina1969method} and matrix Krasulina updates \citep{tang2019matrixkrasulina}. 
However, the resulting gradient noise typically limits these methods to \emph{sublinear} rates of convergence \citep{balsubramani2013streaming,jain2016streaming,allen2017first,li2018near,huang2021streaming,liang2023optimality}. Variance-reduced eigenspace methods such as VR-PCA and stochastic power variants recover linear convergence \citep{shamir2015stochastic,garber2016faster,xu2018accelerated,kim2020stochastic,shamir2016fast}, but require periodic full-pass computations and therefore will suffer from the same straggler issue again.

The above results leave open whether an eigenspace computation method can preserve linear convergence based on component-wise updates without relying on worker synchronization.
The \emph{Incremental Aggregated Gradient} (IAG) methodology provides a natural mechanism for reconciling component-wise asynchronous updates with linear convergence. At each iteration, IAG refreshes one entry in a persistent gradient table and uses the aggregate of all cached entries as a full-gradient surrogate. It therefore requires only one component evaluation per update while supporting a constant stepsize \citep{blatt2007convergent}. Under bounded staleness $\tau$, components may be refreshed cyclically, in a network-imposed order, or upon worker arrivals, allowing a parameter server to update without waiting for all workers \citep{aytekin2016analysis}.
For smooth strongly convex finite sums, IAG and its proximal extensions achieve linear convergence \citep{gurbuzbalaban2017convergence,vanli2018global}. Extensions of IAG cover nonconvex proximal problems under an error-bound condition \citep{peng2019nonconvex} and relatively smooth problems under Bregman distance growth \citep{zhang2021proximal}. Stochastic variants further allow noisy component-gradient oracles in Euclidean settings \citep{wang2023linear,wang2024dual}. 
Algorithmically, these methods rely on stored component directions sharing a common linear representation and hence being directly aggregable, a property that is not intrinsic to the Grassmannian setting.
To our knowledge, IARG is the closest existing asynchronous algorithm specifically developed for leading-eigenvector computation but is restricted to $k=1$ \citep{wang2023incremental}.
Its iterate is a single unit vector updated by normalization, so the method as proposed neither represents a $k$-dimensional subspace nor enforces block orthogonality. Recovering a leading-$k$ subspace through repeated \emph{deflation} would require $k$ separate runs and a chain of eigengaps $\lambda_r(\bm A)>\lambda_{r+1}(\bm A)$ for all $r\le k$, whereas direct subspace computation must handle right-orthogonal quotient invariance under only the target gap $\lambda_k(\bm A)>\lambda_{k+1}(\bm A)$.

\paragraph{Theoretical gap.}
Transferring the IAG algorithmic principle to eigenspace computation is \emph{not} a routine application of existing works on Euclidean or Bregman IAG.
There are two intertwined challenges.
First, the variable is a subspace on the Grassmannian quotient rather than a vector in a common linear space.
Consequently, stale component Riemannian gradients are computed at past subspaces and therefore live in different tangent spaces; aggregating them without transport requires an extrinsic update whose intrinsic effect must be justified.
Second, the eigenspace objective does not possess a uniform global gradient-dominance constant: the useful error bound becomes uniform only inside a suitable basin around the target eigenspace.
An asynchronous analysis must therefore prove not only descent under delayed information, but also basin invariance for the delayed recursion itself.
These coupled geometric and dynamical issues are not resolved by existing IAG and PLIAG analyses, which rely on component gradients evaluated at different iterates belonging to the \emph{same} vector space and on a growth or error-bound condition that holds \emph{uniformly} over the relevant region. Neither property is automatic in the eigenspace computation problem.
IARG treats stale aggregation on the sphere, but its analysis is tailored to scalar alignment with a leading eigenvector and requires initialization within an angle of $\pi/4$ from it. Expressed in terms of a target angular error $\varepsilon>0$, measured by the angle between the output direction and the leading eigenvector, IARG has iteration complexity $\mathcal O(\tau\delta^{-2}\log(1/\varepsilon))$. This bound depends quadratically on the inverse eigengap $\delta^{-1}$, with hidden constants determined by absolute spectral levels and data-norm bounds. 
Together, these limitations expose an important theoretical gap: a Grassmannian convergence theory for asynchronous leading-$k$ subspace dynamics that offers sharper spectral dependence and constants invariant under spectral shifts remains unavailable.

\subsection{Our Contributions}
We develop {\algname} (\textbf{GRASS}mannian \textbf{I}ncremental \textbf{A}ggregation) for asynchronous finite-sum eigenspace computation. Our main contributions are as follows.

\begin{enumerate}[label=$\bullet$, leftmargin=12pt, topsep=0pt]
	\item \emph{Transport-free aggregation on the Grassmannian.}
		We introduce a stale-gradient aggregation scheme tailored to the quotient geometry of leading-$k$ eigenspace computation.
		{\algname} updates a Stiefel representative by a polar step driven by a stale ambient gradient aggregate, thereby avoiding vector transport between tangent spaces.
		Despite its extrinsic form, we prove that the update's first-order Grassmannian motion is precisely the negative projection of the stale aggregate onto the current tangent space.
		The method also computes the full leading subspace directly, avoiding deflation and requiring only the target eigengap rather than a chain of gaps.
		With overhead proportional to the number of refreshed components, it supports a parameter-server implementation that incorporates asynchronous worker returns without synchronized rounds.

	\item \emph{Tight angle-dependent gradient-dominance geometry.}
		We derive an angle-dependent gradient-dominance (or Riemannian Polyak--{\L}ojasiewicz (PL)) inequality for Problem \eqref{eq:eigen}:
	    \[
		    \frac12\|\grad F(\bm W)\|_F^2
		    \ge
		    \rho(\delta,\theta_k)\,\bigl(F(\bm W)-F^*\bigr),
	    \]
	    where $\theta_k$ denotes the largest principal angle between the current subspace and target eigenspace, with $\sin\theta_k$ being their projection spectral distance.
	    In particular, we show that the best possible coefficient depending on this angle is $\rho_*(\delta,\theta_k)=2\delta\cos^2\theta_k$ under only the eigengap condition $\delta > 0$.
	    Thus, the lack of a global PL constant is quantified exactly: the angle-dependent PL coefficient deteriorates monotonically to zero as $\theta_k$ approaches $\pi/2$, at the sharp rate $\cos^2\theta_k$.
	    This result is of independent interest as it provides a standalone error-bound for analyzing first-order methods for leading-eigenspace computation and avoids unnecessary slack in basin and rate estimates.

	\item \emph{Two-phase linear convergence via basin invariance.}
	    We prove that, despite stale component information, the iterates remain in a basin where the sharp angle-dependent gradient-dominance bound stays uniformly effective.
	    This yields a two-phase linear convergence guarantee under an eigengap $\delta>0$ and maximum staleness $\tau$. 
		For a target accuracy $\varepsilon>0$,\footnote{Accuracy is measured by the squared Grassmannian geodesic distance between the output subspace and the leading-$k$ eigenspace; see \eqref{eq:grass-angle}.} the iteration complexity is $\mathcal O(\tau\delta^{-3/2}\log(1/(\delta\varepsilon)))$ for initialization in a prescribed broad basin, and improves to $\mathcal O((\tau+1)\delta^{-1}\log(1/\varepsilon))$ once the iterates enter the local region. 
		The constants depend on component spectral spreads rather than absolute spectral levels, making the guarantees invariant under shifts $\bm A_i\mapsto\bm A_i+c_i\bm I$ for arbitrary scalars $c_i\in\mathbb R$.
\end{enumerate}

\begin{table}[!t]
	\centering
	\scriptsize
	\setlength{\tabcolsep}{2.0pt}
	\renewcommand{\arraystretch}{2.1}
	\begin{minipage}{0.95\textwidth}
		\caption{Comparison of related methods. \emph{Setting} specifies the problem domain or target of each method; \emph{Incremental} and \emph{Asynchronous} indicate whether a method uses component-wise incremental updates and supports asynchronous parameter-server updates. One oracle call computes one component gradient. For {\algname}, the reported asymptotic high-accuracy oracle complexity corresponds to the implementation where each iteration makes one oracle call to refresh one component gradient; staleness factors are written for the regime $\tau\ge1$, while the formal theorem covers $\tau=0$. Constants such as $k$-dependent factors, spectral-scale factors, and component-norm bounds are suppressed.}
		\label{tab:method-comparison}
		\vspace{0.1em}
		\begin{tabular}{@{}
			>{\raggedright\arraybackslash}p{0.26\linewidth}
			>{\centering\arraybackslash\hspace*{-0.3\linewidth}}m{0.17\linewidth}
			>{\centering\arraybackslash\hspace*{-0.2\linewidth}}m{0.16\linewidth}
			>{\centering\arraybackslash\hspace*{-0.1\linewidth}}m{0.12\linewidth}
			>{\raggedright\arraybackslash\hspace*{0.2\linewidth}}m{0.25\linewidth}
			@{}}
			\toprule
			Method & \emph{Setting} & \emph{Incremental} & \emph{Asynchronous} & \emph{Oracle Complexity} \\
			\midrule
			Orthogonal Iteration \citep{saad2011numerical} & $k$-dim eigenspace & \xmark & \xmark & $\mathcal O\!\left(n\delta^{-1}\log(1/\varepsilon)\right)$\\
			\midrule
			RGD \citep{alimisis2024geodesic} & $k$-dim eigenspace & \xmark & \xmark & $\mathcal O\!\left(n\delta^{-1}\log(1/\varepsilon)\right)$ \\
			\midrule
			Block Oja \citep{huang2021streaming} & $k$-dim eigenspace & \cmark & \xmark & $\mathcal O\!\left(\delta^{-2}/\varepsilon\right)$ \\
			\midrule
			VR-PCA \citep{shamir2016fast} & $k$-dim eigenspace & \cmark & \xmark & $\mathcal O\!\left((n+\delta^{-2})\log(1/\varepsilon)\right)$ \\
			\midrule
			FedPower \citep{guo2024fedpower} & $k$-dim eigenspace & \xmark & \xmark & $\mathcal O\!\left(n\delta^{-1}\log(1/\varepsilon)\right)$ \\
			\midrule
			IARG \citep{wang2023incremental} & $1$-dim eigenvector & \cmark & \cmark & $\mathcal O\!\left(\tau\delta^{-2}\log(1/\varepsilon)\right)$\\
			\midrule
			\rowcolor{gray!10}
			{\algname} (this paper) & $k$-dim eigenspace & \cmark & \cmark & $\mathcal O\!\left(\tau\delta^{-1}\log(1/\varepsilon)\right)$\\
			\bottomrule
		\end{tabular}
	\end{minipage}
\end{table}

\paragraph{Comparison with related work.}
Table~\ref{tab:method-comparison} compares our contribution to existing works for eigenspace or eigenvector computation under a common component-oracle accounting.
Full-operator or round-based methods including Orthogonal Iteration, RGD, and FedPower achieve linear eigengap-dependent rates, but they are synchronous methods and do not update from individual component arrivals.
Block Oja provides cheap incremental updates but only sublinear convergence, while VR-PCA gives a variance-reduced linear rate but does not support stale component arrivals and retains a $\delta^{-2}$ eigengap term.
IARG treats only the $k=1$ case; its available iteration bound also has the looser $\delta^{-2}$ eigengap dependence and constants tied to absolute spectral levels, hence is not shift-invariant.
Beyond the direct eigenspace comparisons in Table~\ref{tab:method-comparison}, IAG and PLIAG provide methodological benchmarks for stale incremental aggregation. For smooth and strongly convex Euclidean problems, IAG attains an iteration complexity of $\mathcal O\!\left(\tau \kappa \log(1/\varepsilon)\right)$, governed by the condition number $\kappa$ \citep{vanli2018global}. In the Grassmannian eigenspace problem, the eigengap $\delta$ plays the corresponding curvature role, so the local $\delta^{-1}$ dependence of {\algname} parallels IAG's linear dependence on $\kappa$.
PLIAG provides the corresponding Bregman benchmark: its iteration complexity of $\mathcal O\!\left(\tau\phi(\tau)\kappa\log(1/\varepsilon)\right)$ is linear in $\kappa$, where $\phi(\tau)$ is a monotonically increasing function defined by the Bregman-distance bound in \cite{zhang2021proximal}.
Taken together, these comparisons highlight that {\algname} directly computes the full leading-$k$ eigenspace while providing, to our knowledge, the first linear convergence theory for asynchronous optimization on the Grassmannian.

\subsection{Notation and Preliminaries}\label{sec:notation}

\paragraph{Basic notation.}
For $m\in\mathbb{N}$, write $[m]\coloneqq\{1,\ldots,m\}$.
Let $\langle \bm X,\bm Y\rangle\coloneqq\tr(\bm X^\top\bm Y)$, and use $\|\cdot\|_2$ and $\|\cdot\|_F$ for spectral and Frobenius norms.
Let $\bm W^*$ be any global minimizer of \eqref{eq:eigen}, and let $F^*$ be the optimal value.
For $\bm B\in\mathbb S^d$, write
\(
	\lambda_1(\bm B)\ge\lambda_2(\bm B)\ge\cdots\ge\lambda_d(\bm B)
\)
for its eigenvalues in nonincreasing order.
Let $\bm A=\bm U\bm\Lambda\bm U^\top$, where $\bm\Lambda=\diag(\lambda_1,\ldots,\lambda_d)$ and $\lambda_1\ge\cdots\ge\lambda_d$.
Set
\[
	\bm{\Lambda}_k \coloneqq \diag(\lambda_1,\ldots,\lambda_k),
	\qquad
	\bar{\bm{\Lambda}}_{d-k} \coloneqq \diag(\lambda_{k+1},\ldots,\lambda_d).
\]
Define
\[
	\delta \coloneqq \lambda_k(\bm A)-\lambda_{k+1}(\bm A),
	\qquad
	\nu \coloneqq \lambda_1(\bm A)-\lambda_d(\bm A).
\]
For each component matrix, set
\[
	\nu_i \coloneqq \lambda_1(\bm A_i)-\lambda_d(\bm A_i),
	\qquad
	\nu_{\rm avg} \coloneqq \frac{1}{n}\sum_{i=1}^n \nu_i.
\]
For any full-column-rank $\bm X\in\mathbb R^{d\times k}$, its unique thin polar decomposition is given by \citep{higham1986computing}
\[
	\bm X=\Polar(\bm X)(\bm X^\top\bm X)^{1/2},
	\qquad
	\Polar(\bm X)\coloneqq \bm X(\bm X^\top\bm X)^{-1/2}\in\St(d,k),
\]
where $\Polar(\bm X)$ is the orthonormal polar factor.

\paragraph{Grassmann geometry.}
Denote the Grassmannian, i.e., the set of $k$-dimensional subspaces of $\mathbb R^d$, by $\Gr(d,k)$, which is identified with the quotient space $\St(d,k)/\mathbb O(k)$.
For $\bm W\in\St(d,k)$, write $[\bm W]\coloneqq\{\bm W\bm Q:\bm Q\in\mathbb O(k)\} \in \Gr(d,k)$.
For $\bm W,\bm W'\in\St(d,k)$, let 
\(
	0\le\theta_1(\bm W,\bm W')\le\cdots\le\theta_k(\bm W,\bm W')\le\pi/2
\)
be the principal angles between $[\bm W]$ and $[\bm W']$.
Then, the Grassmannian geodesic distance between $[\bm W]$ and $[\bm W']$ is
\begin{equation} \label{eq:grass-angle}
	\dist_{\Gr}(\bm W,\bm W')
	\coloneqq
	\left(\sum_{\ell=1}^k \theta_\ell^2(\bm W,\bm W')\right)^{1/2}.
\end{equation}
The tangent space at $[\bm W]$ is
\(
	\mathcal T_{[\bm W]}\Gr(d,k)
	\coloneqq
	\{\bm\Xi\in\mathbb R^{d\times k}:\bm W^\top\bm\Xi=\bm 0\},
\)
and the tangent projection of $\bm G\in\mathbb R^{d\times k}$ is
\[
	\mathscr P_{\mathcal T_{[\bm W]}\Gr(d,k)}(\bm G)
	\coloneqq
	(\bm I-\bm W\bm W^\top)\bm G.
\]
% \paragraph{Polar factor and logarithm map}
For $\bm W,\bm W'\in\St(d,k)$ with $\bm W^\top\bm W'$ invertible, set
\(
\bm L\coloneqq \bm W'(\bm W^\top\bm W')^{-\!1}-\bm W
\)
and let $\bm L=\bm R\bm\Sigma\bm V^\top$ be a thin SVD with $\bm R \in \St(d,k)$ and $\bm V \in \mathbb{O}(k)$.
The Grassmann logarithm \citep{absil2004riemannian} is
\[
	\Log_{\bm W}(\bm W')\coloneqq
	\bm R\,\arctan(\bm\Sigma)\,\bm V^\top,
\]
where $\arctan(\bm\Sigma)$ is applied entrywise.
Then $\Log_{\bm W}(\bm W')$ lies in $\mathcal T_{[\bm W]}\Gr(d,k)$ and has Frobenius norm $\dist_{\Gr}(\bm W,\bm W')$.

\section{The {\algname} Method}

We consider finite-sum objectives represented over the Stiefel manifold,
$F(\bm W)=\frac{1}{n}\sum_{i=1}^n F_i(\bm W)$ with $\bm W\in\St(d,k)$, where the underlying decision variable is the subspace $[\bm W]\in\Gr(d,k)$.
For the leading eigenspace computation problem, the component objective is $F_i(\bm W)=-\tr(\bm W^\top \bm A_i \bm W)$ and the corresponding Riemannian gradient over $\Gr(d,k)$ \citep{absil2008optimization,hu2020brief,boumal2023introduction} admits the closed form
\begin{align}
	\grad F_i(\bm W)
	=
	\mathscr{P}_{\mathcal{T}_{[\bm W]}\Gr(d,k)}\left(\nabla F_i(\bm W)\right)
	=
	-2\left(\bm I-\bm W\bm W^\top\right)\bm A_i \bm W,
	\label{eq:local-riem-grad}
\end{align}
where $\nabla F_i(\bm W)=-2\bm A_i\bm W$ is the Euclidean gradient of $F_i$ at $\bm W$.

\subsection{Design Principles}

The construction of {\algname} is guided by two requirements: preserving the computational economy of incremental aggregation and respecting the geometry of eigenspace computation.
The first requires that, when only a small subset of components is refreshed, the aggregate can be corrected using only the newly arrived information.
The geometric requirement complicates this recursion: a cached component gradient computed at an earlier iterate $\bm W^s$ belongs to $\mathcal T_{[\bm W^s]}\Gr(d,k)$, whereas the update is formed at the current subspace $[\bm W^t]$.
Because cached gradients may have different base points, they are \emph{not} intrinsically summable.
Transporting each of them to the current tangent space by parallel transport is geometrically natural, but the transport destination \emph{changes} with every iterate.
Consequently, even unrefreshed cache entries must be represented anew, and a direct implementation would require transporting all $n$ cached gradients at every iteration, with each Grassmannian transport costing $\mathcal O(dk^2+k^3)$; this would forfeit the cost advantage of component-wise updates.

Rather than transporting the cached gradients to a common tangent space, {\algname} treats their matrix representatives as ambient quantities, aggregates them incrementally in $\mathbb R^{d\times k}$, and applies the polar map to the resulting ambient step.
Caching Riemannian rather than Euclidean component gradients is deliberate: it removes the normal component at the subspace where each gradient is computed.
Because the normal space varies with the base point, retaining such a component could induce a spurious first-order change in the current subspace when stale information is used at a later iterate.
Caching Riemannian gradients removes such past normal components before aggregation.
The following construction implements these design choices.

\subsection{Algorithmic Construction}
Let $\bm G_i^t$ denote the cached component gradient for index $i$ that is available at the beginning of iteration $t$.
We initialize the table at the initial point by setting
\(
	\bm G_i^{0} \coloneqq \grad F_i(\bm W^0)
\)
for all $i\in[n]$ and $\bm G^{0}=\frac1n\sum_{i=1}^n \bm G_i^{0}$.
We write $\tau_i(t)\in\{0,\ldots,t\}$ for the staleness of the cached entry
\(
	\bm G_i^t \coloneqq \grad F_i(\bm W^{t-\tau_i(t)}).
\)
Initially, $\tau_i(0)=0$ for all $i$.
At iteration $t=0,1,2,\ldots$, the aggregated direction is formed by averaging the current table:
\begin{align}
	\bm G^{t}
	\coloneqq
	\frac{1}{n}\sum_{i=1}^{n}\bm G_i^{t}
	= \frac{1}{n}\sum_{i=1}^{n}\grad F_i(\bm W^{t-\tau_i(t)}).
	\label{eq:approx-rg-gen}
\end{align}
Although $\bm G^t$ is not an intrinsic tangent vector at $[\bm W^t]$, its ambient matrix representation defines the polar update
\begin{align}
	\bm W^{t+1}
	=
	\Polar(\bm W^{t}-\eta \bm G^{t}),
	\quad t=0,1,2,\ldots.
	\label{eq:iter}
\end{align}
Here, $\Polar(\cdot)$ maps a full-column-rank matrix to its polar factor, i.e., the closest point on $\St(d,k)$ under the Frobenius norm.
Thus, the algorithm takes an ambient-space step and restores orthonormality without transporting stale gradients across tangent spaces.
This projection guarantees feasibility, but its effect on the Grassmannian is not immediate because $\bm G^t$ is not tangent at $[\bm W^t]$.
Lemma~\ref{lem:log-decomp} in the next section provides the missing geometric link: the induced first-order motion is governed by the tangent projection $\mathscr{P}_{\mathcal{T}_{[\bm W^t]}\Gr(d,k)}(\bm G^t)$, with only a higher-order discrepancy.

After $\bm W^{t+1}$ is computed, a nonempty subset $\mathcal S_{t+1}\subseteq[n]$ of component indices is selected or arrives with newly available gradients
\begin{align}
	\widehat{\bm G}_i^{t+1}
	= \grad F_i (\bm W^{t+1-\tau_i(t+1)}),
	\qquad i\in \mathcal S_{t+1}.
\end{align}
The table is then updated by setting $\bm G_i^{t+1}=\widehat{\bm G}_i^{t+1}$ for $i\in\mathcal S_{t+1}$ and $\bm G_i^{t+1}=\bm G_i^t$ otherwise.
For each refreshed entry, $\tau_i(t+1)$ records the age of the model at which the newly arrived component gradient was computed.
For the unrefreshed entries, the same cached gradients are retained, which is equivalently represented by increasing their staleness by one:
\(
	\tau_i(t+1)=\tau_i(t)+1
	\;\text{for}\; i\notin \mathcal S_{t+1}.
\)
The singleton choice $|\mathcal S_{t+1}|=1$ gives the fully incremental implementation, whereas the choice $\mathcal S_{t+1}=[n]$ with $\tau_i(t+1)=0$ for all $i$ recovers the synchronous specialization.
Because all unrefreshed entries remain unchanged, the aggregate itself admits the incremental recursion
\begin{align}
	\bm G^{t+1}
	= \bm G^{t}
	+\frac{1}{n}\sum_{i\in \mathcal S_{t+1}}
	\left(\widehat{\bm G}_i^{t+1}-\bm G_i^t\right),
	\quad t=0,1,2,\ldots. 
	\label{eq:agg-recursion}
\end{align}

The unified algorithmic procedure of {\algname} for eigenspace computation is summarized in Algorithm~\ref{alg:generic}.
It encompasses serial incremental, asynchronous distributed, and synchronous distributed implementations through different choices of the refreshed subset \(\mathcal S_{t+1}\) and the corresponding stalenesses $\{\tau_i(t+1):\;i\in \mathcal S_{t+1}\}$.

\begin{algorithm}[t]
	\caption{{\algname} for Eigenspace Computation}
	\label{alg:generic}
	\begin{algorithmic}[1]
		\STATE \textbf{Input:} stepsize $\eta>0$, initial point $\bm W^{0}\in \St(d,k)$
		\STATE \textbf{Initialization:} for each $i\in[n]$, set $\bm G_i^{0}=\grad F_i(\bm W^{0})$; set $\bm G^{0}=\frac{1}{n}\sum_{i=1}^{n}\bm G_i^{0}$
		\FOR{$t=0,1,2,\dots$}
		\STATE $\bm W^{t+1}=\Polar\left(\bm W^{t}-\eta \bm G^{t}\right)$
		\STATE Receive/select a nonempty subset $\mathcal S_{t+1}\subseteq[n]$ with stalenesses $\{\tau_i(t+1):i\in \mathcal S_{t+1}\}$
		\STATE Compute the refreshed gradients
		\(
			\widehat{\bm G}_i^{t+1} = \grad F_i\left(\bm W^{t+1-\tau_i(t+1)}\right) \text{ for } i\in \mathcal S_{t+1}
		\)
		\STATE Update the aggregated direction:
		\(
			\bm G^{t+1} = \bm G^{t} + \frac{1}{n}\sum_{i\in \mathcal S_{t+1}}(\widehat{\bm G}_i^{t+1} - \bm G_i^{t})
		\)
		\STATE Set $\bm G_i^{t+1}=\widehat{\bm G}_i^{t+1}$ for all $i\in \mathcal S_{t+1}$ and $\bm G_i^{t+1}=\bm G_i^{t}$ for all $i\notin \mathcal S_{t+1}$
		\ENDFOR
	\end{algorithmic}
\end{algorithm}

\paragraph{Asynchronous implementation.}
In a distributed implementation, worker $i$ owns $F_i$ and its local cache $\bm G_i^t$, while the server stores only $\bm W^t$ and the aggregate $\bm G^t$; hence the full $\mathcal O(ndk)$ table need not be stored centrally.
Whenever a subset $\mathcal S_{t+1}$ returns, each worker sends only the correction $(\widehat{\bm G}_i^{t+1}-\bm G_i^t)/n$, and the server applies \eqref{eq:agg-recursion} immediately without waiting for the remaining workers.
Faster workers may therefore contribute more frequently, while information from slower workers is incorporated upon arrival.
For $|\mathcal S_{t+1}|$ returned components, updating the aggregate costs $\mathcal O(|\mathcal S_{t+1}|dk)$, and the polar update in \eqref{eq:iter} costs $\mathcal O(dk^2+k^3)$ arithmetic operations by forming the $k\times k$ Gram matrix, computing its inverse square root, and multiplying back, so the server-side update overhead remains \emph{independent} of $n$ apart from the arriving corrections.

\section{Theoretical Analysis}\label{sec:main-results}

This section presents the main theoretical guarantees for {\algname}. The detailed proofs and the supporting technical results are deferred to the Appendix.

\paragraph{Roadmap of the analysis.}
The proof is organized around the fact that the eigenspace objective has useful curvature only when the trajectory is in the right angular basin.
We first establish exact projector identities for the gradient energy and value gap, which yield a tight angle-dependent gradient-dominance inequality.
The linear rate is therefore not a direct consequence of a generic gradient-dominance theorem: there is no uniform global gradient-dominance constant, and the local constant deteriorates as the largest principal angle approaches $\pi/2$.
Thus, the convergence proof must simultaneously show descent and prove that the delayed iterates never leave a region where this geometry is effective.
A second ingredient concerns the algorithmic structure.
{\algname} does not move along the current Riemannian gradient; it uses an ambient aggregate of component gradients computed at past subspaces, which live in different tangent spaces.
We therefore analyze the extrinsic polar update through its Grassmann logarithm, identify the current tangent component that drives first-order descent, and bound the stale-gradient mismatch by a recursion over recent aggregate-gradient norms.
A Lyapunov basin-invariance argument then absorbs this delay penalty and keeps the trajectory inside an angular neighborhood where the gradient-dominance coefficient is uniformly positive.
With this closed loop in place, the final step solves the delayed recursion, yielding the broad-basin and sharper local rates.

The convergence theory is built around two structural assumptions.  The first is the standard eigengap condition for identifying the leading-$k$ invariant subspace, and the second is the bounded-delay condition used to model stale component information in the incremental aggregation table.

\begin{assumption}\label{as:eigengap}
	The eigengap between the $k$-th and $(k+1)$-th eigenvalues of $\bm A$ is positive, i.e.,
	\(
	\delta\coloneqq \lambda_k-\lambda_{k+1} > 0.
	\)
\end{assumption}

\begin{assumption}\label{as:delay}
	There exists a constant $\tau \in \mathbb{N}$ such that $\tau_i(t)\le \tau$ for all $i\in[n]$ and all $t$.
\end{assumption}
This assumption covers the serial incremental and asynchronous distributed implementations, where stale components generally lead to $\tau \ge 1$ when $n>1$, as well as the synchronous distributed specialization, where all component gradients are refreshed at the current iterate and hence $\tau=0$.

\subsection{Eigenspace Geometry}\label{sec:geometry}

The next two lemmas quantify the geometry of the trace objective.
They relate the objective gap and the gradient norm to the principal angles between the current subspace and the target eigenspace.
Together, these identities give a tight angle-dependent gradient-dominance bound.
For notational convenience, we first define $\bm{Z}\coloneqq \bm{U}^\top \bm{W}=\begin{bmatrix}\bm{Z}_1\\ \bm{Z}_2\end{bmatrix}$ with $\bm{Z}_1\in\mathbb{R}^{k\times k}$ and $\bm{Z}_2\in\mathbb{R}^{(d-k)\times k}$ for $\bm{W}\in\St(d,k)$.

\begin{lemma}[Gradient-energy identity]\label{lem:grad-energy}
	Let $\bm{W}\in\St(d,k)$ and $\bm{P}\coloneqq \bm{Z}\bm{Z}^\top$ be the rank-$k$ projector. Block $\bm{P}$ as
	\begin{align}
		\bm{P}=\begin{bmatrix}\bm{P}_{11}&\bm{P}_{12}\\ \bm{P}_{21}&\bm{P}_{22}\end{bmatrix}
		\;\text{with}
		\;\bm{P}_{11}\coloneqq\bm{Z}_1\bm{Z}_1^\top,
		\;\bm{P}_{22}\coloneqq\bm{Z}_2\bm{Z}_2^\top,
		\;\bm{P}_{12}\coloneqq\bm{Z}_1\bm{Z}_2^\top,
		\;\bm{P}_{21}\coloneqq\bm{Z}_2\bm{Z}_1^\top.
		\nonumber
	\end{align}
	Then, the following identity holds:
	\begin{align}
		\frac{1}{4} \left\|\grad F(\bm{W})\right\|_F^2
		=\frac12 \left\|\comm{\bm{\Lambda}_k}{\bm{P}_{11}}\right\|_F^2
		+\frac12 \left\|\comm{\bar{\bm{\Lambda}}_{d-k}}{\bm{P}_{22}}\right\|_F^2
		+ \left\|\bm{\Lambda}_k\bm{P}_{12}\!-\!\bm{P}_{12}\bar{\bm{\Lambda}}_{d-k}\right\|_F^2,
		\nonumber
	\end{align}
	where $\comm{\bm{X}}{\bm{Y}} \coloneqq \bm{X}\bm{Y} - \bm{Y}\bm{X}$ denotes the commutator of square matrices $\bm{X}$ and $\bm{Y}$.
\end{lemma}

\begin{lemma}[Value-gap identity]\label{lem:value-gap}
	Let $\bm{W}\in\St(d,k)$ and $\theta_1,\dots,\theta_k$ be the principal angles between $[\bm{W}]$ and $[\bm{W}^*]$.
	Consider the following two cases of the CS decomposition of $\bm{Z}$ according to the relative size of blocks $\bm{Z}_1$ and $\bm{Z}_2$:

	i) If $d-k\ge k$, then write the CS decomposition of $\bm{Z}$ as
	\begin{align}
		\bm{Z}_1=\bm{Q}_1 \bm{C} \bm{R}^\top,\qquad
		\bm{Z}_2 =
		\bm{Q}_2
		\begin{bmatrix}
			\bm S \\
			\bm 0_{(d-2k)\times k}
		\end{bmatrix}
		\bm{R}^\top,
		\nonumber
	\end{align}
	where $\bm{Q}_1\in\mathbb{O}(k)$, $\bm{Q}_2\in\mathbb{O}(d-k)$, $\bm{R}\in\mathbb{O}(k)$, and
	\begin{align}
		\bm{C}=\diag(\cos\theta_1,\dots,\cos\theta_k),
		\qquad
		\bm{S}=\diag(\sin\theta_1,\dots,\sin\theta_k).
		\nonumber
	\end{align}

	ii) If $d-k<k$ (the first $2k-d$ principal angles are zero),  then write the CS decomposition of $\bm{Z}$ as
	\begin{align}
		\bm{Z}_1
		 & =\bm{Q}_1
		\begin{bmatrix}
			\bm{I}_{2k-d} & \bm{0}   \\
			\bm{0}        & \bm{C}_0
		\end{bmatrix}
		\bm{R}^\top,
		\qquad
		\bm{Z}_2
		=\bm{Q}_2
		\begin{bmatrix}
			\bm{0}_{(d-k)\times(2k-d)},\; \bm{S}_0
		\end{bmatrix}
		\bm{R}^\top,
		\nonumber
	\end{align}
	where $\bm{Q}_1\in\mathbb{O}(k)$, $\bm{Q}_2\in\mathbb{O}(d-k)$, $\bm{R}\in\mathbb{O}(k)$,
	\(
	\bm{C}_0=\diag(\cos\theta_{2k-d+1},\dots,\cos\theta_k),
	\) 
	and
	\(
	\bm{S}_0=\diag(\sin\theta_{2k-d+1},\dots,\sin\theta_k).
	\)
	In this case, we write
	\begin{align}
		\bm{C}
		=
		\begin{bmatrix}
			\bm{I}_{2k-d} & \bm{0}   \\
			\bm{0}        & \bm{C}_0
		\end{bmatrix}
		\in \mathbb{R}^{k\times k},
		\qquad
		\bm{S}
		=
		\begin{bmatrix}
			\bm{0}_{2k-d} & \bm{0}    \\
			\bm{0}        & \bm{S}_0,
		\end{bmatrix}
		\in \mathbb{R}^{k\times k}.
		\nonumber
	\end{align}
	Define two $k\times k$ matrices:
	\begin{align}
		\bm{M}_k \coloneqq \bm{Q}_1^\top \bm{\Lambda}_k \bm{Q}_1,
		\;\;
		\overline{\bm{M}}_k
		\coloneqq
		\begin{cases}
			\begin{bmatrix}
				\bm I_k & \bm 0_{k\times(d-2k)}
			\end{bmatrix}
			\bm{Q}_2^\top \bar{\bm{\Lambda}}_{d-k} \bm{Q}_2
			\begin{bmatrix}
				\bm I_k \\
				\bm 0_{(d-2k)\times k}
			\end{bmatrix},
			 & d-k\ge k, \\[2mm]
			\begin{bmatrix}
				\lambda_{k+1}\bm{I}_{2k-d} & \bm{0}                                          \\
				\bm{0}                     & \bm{Q}_2^\top \bar{\bm{\Lambda}}_{d-k} \bm{Q}_2
			\end{bmatrix},
			 & d-k<k.
		\end{cases}
		\nonumber
	\end{align}
	Then, one has the exact identity
	\begin{align}
		F(\bm{W})-F^*
		 & =\tr\left((\bm{M}_k-\overline{\bm{M}}_k) \bm{S}^2\right).
		\nonumber
	\end{align}
\end{lemma}

Lemmas~\ref{lem:grad-energy} and \ref{lem:value-gap} isolate the two quantities that drive the local geometry: the gradient energy is controlled by spectral commutators, while the objective gap is controlled by the squared sine factors from the CS decomposition.
As a consequence, the value-gap identity implies the usual eigengap error bound between objective gap and Grassmannian geodesic distance.

\begin{proposition}[Tight angle-dependent gradient dominance]\label{lem:localPL}
	Suppose that Assumption~\ref{as:eigengap} holds. Then, for every $\bm W\in\St(d,k)$, the Riemannian gradient satisfies
	\begin{align}
		\frac{1}{2} \|\grad F(\bm{W})\|_F^2 \ge 2 \delta \cos^2 (\theta_k) \left(F(\bm{W})-F^*\right),
		\label{eq:localPL}
	\end{align}
	where $\theta_k=\theta_k(\bm W,\bm W^*)$ is the largest principal angle between the two subspaces.
	Moreover, the factor $2\delta\cos^2\theta_k$ is pointwise tight; in specific, it cannot be replaced by any uniformly larger pointwise coefficient depending only on $\delta$ and  $\theta_k$.
\end{proposition}

Proposition~\ref{lem:localPL} shows that Problem \eqref{eq:eigen} is naturally governed by an \emph{angle-dependent gradient-dominance} property. Indeed, the factor $2\delta \cos^2\theta_k$ depends explicitly on the largest principal angle between the current subspace $[\bm W]$ and the target eigenspace $[\bm W^*]$.
Since $\|\bm W\bm W^\top-\bm W^*(\bm W^*)^\top\|_2=\sin\theta_k$, the largest principal angle $\theta_k$ is the angular representation of the projection spectral distance between the two subspaces.
The tightness statement shows that this angle dependence is not merely a proof artifact, but reflects an intrinsic geometric feature of the problem under only the eigengap condition. 
The tight coefficient also clarifies the relation to the result by \citet[Proposition 4]{alimisis2024geodesic}, which gives a related gradient-dominance coefficient proportional to
\(
c_Q(\theta_k/\tan\theta_k)^2
\)
with $c_Q=4/\pi^2$.
Proposition~\ref{lem:localPL} instead yields the tight coefficient $\cos^2\theta_k$, which is pointwise no smaller because
\(
	\begin{aligned}
		\frac{4}{\pi^2}\left(\frac{\theta_k}{\tan\theta_k}\right)^2
		= \cos^2(\theta_k)\cdot\frac{4\theta_k^2}{\pi^2\sin^2\theta_k}
		\le \cos^2\theta_k.
	\end{aligned}
\)
The inequality uses $\sin x\ge 2x/\pi$ for $x\in[0,\pi/2]$.
The sharper coefficient is obtained by exploiting the eigenspace-specific gradient-energy and value-gap identities in our proof.

In general, Problem~\eqref{eq:eigen} does \emph{not} satisfy a global gradient-dominance inequality, i.e. there does not exist an absolute constant $\mu>0$ such that
\(
	\frac{1}{2}\|\grad F(\bm W)\|_F^2 \ge \mu\big(F(\bm W)-F^*\big)
	\text{ for all } \bm W\in\St(d,k).
\)
The reason is that any $k$-dimensional invariant subspace of $\bm A$ is a stationary point of $F$, since $\grad F(\bm W)= -2(\bm I-\bm W\bm W^\top)\bm A\bm W=\bm 0$ whenever $\bm A\,\mathrm{span}(\bm W)\subseteq \mathrm{span}(\bm W)$. If such a subspace is not the leading-$k$ eigenspace, then it is stationary but suboptimal, i.e., $F(\bm W)-F^*>0$. Hence, a global gradient-dominance inequality with a uniform positive constant cannot hold.
Nevertheless, Proposition~\ref{lem:localPL} recovers a uniform PL inequality on any angular region where $\theta_k$ is bounded away from $\pi/2$.
For any $\alpha\in[0,\pi/2)$, let
\(
	\mathcal V_\alpha
	\coloneqq
	\left\{\bm W\in\St(d,k):\theta_k(\bm W,\bm W^*)\le\alpha\right\}.
\)
Since $\theta_k\le\alpha$ on $\mathcal V_\alpha$, the sharp coefficient satisfies $2\delta\cos^2\theta_k\ge 2\delta\cos^2\alpha>0$. Hence, for all $\bm W\in\mathcal V_\alpha$,
\[
	\frac{1}{2}\|\grad F(\bm W)\|_F^2 \ge 2\delta \cos^2(\alpha) \left(F(\bm W)-F^*\right).
\]
Thus, once the iterate stays in such a fixed angular neighborhood of the target eigenspace, the angle-dependent gradient-dominance property reduces to a PL-type inequality with the explicit constant $2\delta\cos^2\alpha$. 
Nevertheless, the bound above does not by itself imply convergence of {\algname}: the update uses a stale ambient aggregate, and the coefficient is useful only while the trajectory stays in $\mathcal V_\alpha$. The next two subsections address these issues through a delayed descent recursion for the actual polar step and a basin-invariance argument.

\subsection{One-Step Analysis}\label{sec:one-step}

Although Algorithm~\ref{alg:generic} is written as a polar retraction step in the ambient matrix space, the convergence argument is carried out on $\Gr(d,k)$.
The following logarithmic expansion shows that, up to a second-order error, the polar step follows the tangent projection of the aggregated direction.
This allows us to use a Riemannian Taylor expansion of $F$ while keeping the constants in terms of the shift-invariant spreads $\nu$ and $\nu_{\rm avg}$.

\begin{lemma}[Decomposition of logarithm map]\label{lem:log-decomp}
	Let $\bm W^t\in\St(d,k)$ and let $\bm G^t\in\mathbb R^{d\times k}$ be the aggregate defined in \eqref{eq:approx-rg-gen}.
	If the step size $\eta > 0$ satisfies
	\begin{align}
		\eta \le 1/(2 \nu_{\rm avg}),
		\label{eq:eta-cond}
	\end{align}
	then the polar step $\bm W^{t+1} \coloneqq \Polar(\bm W^t-\eta \bm G^t)$ is uniquely defined.
	Moreover, the Grassmann logarithm $\Log_{\bm W^t}(\bm W^{t+1})$ is uniquely defined and admits the decomposition
	\begin{align}
		\Log_{\bm W^t}(\bm W^{t+1}) = -\eta \mathscr{P}_{\mathcal{T}_{[\bm W^t]}\Gr(d,k)}(\bm G^t) + \bm E^t
		\quad\text{with}\quad
		\|\bm E^t\|_F \le C_{\rm po} \eta^2 \|\bm G^t\|_F^2,
		\label{eq:log-decomp}
	\end{align}
	where $C_{\rm po} \coloneqq 2+{4\sqrt{k}}/{3}$.
\end{lemma}

The one-step descent estimate below is the main algorithmic inequality.

\begin{proposition}[One-step descent of function value]\label{prop:one-step-descent}
	Let $\bm{W}^t$ and $\bm{W}^{t+1}$ with $t\ge0$ be any two consecutive iterates produced by Algorithm~\ref{alg:generic}. Suppose that Assumption~\ref{as:delay} is satisfied. Then, for any stepsize $\eta>0$ obeying
	\[
		\eta \le \eta_0 \coloneqq \min \left\{ \frac{1}{2\nu_{\rm avg}},\; \frac{1}{8 \nu (C_{\rm po} \sqrt{k} + 4)} \right\},
	\]
	it holds for every $t\ge0$ that
	\[
		F(\bm{W}^{t+1}) - F(\bm{W}^{t})
		\leq - \frac{1}{4} \eta \| \grad F (\bm{W}^{t}) \|_F^2
		- \frac{1}{2} \eta \| \bm{G}^t \|_F^2
		+ 144\nu_{\rm avg}^2\tau\eta^3
		\sum_{s=[t-\tau]_+}^{t-1}\!\| \bm{G}^s \|_F^2.
	\]
\end{proposition}

Proposition~\ref{prop:one-step-descent} is the point at which asynchrony enters the analysis.
It has the usual descent term in the true gradient, an additional dissipative term in the aggregated direction, and a delay penalty that accumulates the stale-gradient error over the most recent $\tau$ updates.
The last term vanishes in the (synchronous) full-gradient case $\tau=0$, where all component gradients are refreshed at the current iterate.

\subsection{Basin Invariance}\label{sec:basin-invariance}

The angle-dependent gradient-dominance inequality from Proposition~\ref{lem:localPL} becomes a uniform PL inequality only while the iterates stay inside a fixed angular neighborhood $\mathcal V_\alpha$.
Thus, before proving linear convergence, we must show that the delayed aggregated iterates do not leave this neighborhood.
Building on the one-step descent estimate in Proposition~\ref{prop:one-step-descent}, the following proposition provides this invariance statement through a Lyapunov function that balances objective decrease against possible movement in Grassmannian geodesic distance.

\begin{proposition}[Basin invariance of {\algname} iterates on $\Gr(d,k)$]\label{prop:max-angle}
	Fix an initial radius $\zeta\in[0,\pi/2)$ and choose an angular radius $\alpha\in(\zeta,\pi/2)$. 
	Let $B_{k,\alpha}\coloneqq {17}/{8}+2C_{\rm po}\sqrt{k} \alpha$, $A_{\delta,\alpha}\coloneqq {9\pi^2\tan\alpha}/(\delta\alpha)$, and $R_{\zeta,\alpha}\coloneqq (\alpha^2-\zeta^2)/(\nu\zeta^2)$ for $\zeta>0$.
	Define the compatibility stepsize
	\begin{align}
		\eta_{\rm comp}(\zeta,\alpha)
		\coloneqq
		\begin{cases}
			+\infty, & \zeta=0, \\[1mm]
			\displaystyle \frac{R_{\zeta,\alpha}}{4B_{k,\alpha}}, & \zeta>0,\ \tau=0, \\[3mm]
			\displaystyle
			\min\left\{
			\frac{R_{\zeta,\alpha}}{4B_{k,\alpha}},\,
			\frac{1}{\nu_{\rm avg}\tau}
			\sqrt{\frac{R_{\zeta,\alpha}}
			{4A_{\delta,\alpha}+576R_{\zeta,\alpha}}}
			\right\}, & \zeta>0,\ \tau\ge1.
		\end{cases}
		\label{eq:eta-comp}
	\end{align}
	Suppose that Assumptions~\ref{as:eigengap} and \ref{as:delay} hold, and the stepsize $\eta>0$ satisfy
	\[
		\eta \le \min \left\{
		\begin{aligned}
			\eta_0,\; 
			\frac{1}{4 \nu_{\rm avg} C_{\rm po} \sqrt{k}},\;
			\frac{1}{24 \nu_{\rm avg} (\tau+1)},\,
			\eta_{\rm comp}(\zeta,\alpha)
		\end{aligned}
		\right\}.
	\]
	Let $ \{\bm{W}^t\}_{t \ge 0} \subseteq \St(d,k)$ be the sequence produced by Algorithm~\ref{alg:generic} with initial point $\bm W^0$ such that
	\begin{align}
		\dist_{\Gr}(\bm W^0,\bm W^*) \le \zeta,
		\label{eq:init-zeta}
	\end{align}
	Then, the iterate sequence $\{\bm W^t\}_{t \ge 0}$ satisfies
	\(
		\dist_{\Gr}(\bm W^t,\bm W^*)\le\alpha.
	\)
\end{proposition}

The bound $\eta_{\rm comp}(\zeta,\alpha)$ quantifies the buffer between the initialization radius $\zeta$ and the target basin radius $\alpha$.
Once $\zeta$ and $\alpha$ are chosen, taking $\eta\le\eta_{\rm comp}(\zeta,\alpha)$ ensures that the initial Lyapunov value is small enough to keep the iterates inside $\mathcal V_\alpha$.

\subsection{Two-Phase Linear Convergence}\label{sec:two-phase}

We now combine the three ingredients above.
Proposition~\ref{prop:max-angle} keeps the broad-basin iterates in $\mathcal V_\alpha$, Proposition~\ref{lem:localPL} then supplies a uniform PL constant on that region, and Proposition~\ref{prop:one-step-descent} gives the delayed descent recursion.
The first theorem uses this full mechanism from a broad initialization region whose largest principal angle is bounded away from the boundary $\theta_k=\pi/2$.
The second theorem starts from a stronger objective-gap condition, so the basin is preserved more directly and the admissible stepsize recovers the sharper local eigengap dependence.

\subsubsection{Broad-Basin Linear Convergence}\label{sec:broad}

We first state the convergence guarantee under a general distance initialization inside such a broad basin.
The compatibility stepsize keeps the delayed trajectory inside $\mathcal V_\alpha$, where the gradient dominance holds uniformly.

\begin{theorem}[Broad-basin linear convergence]\label{thm}
	Fix $\zeta\in[0,\pi/2)$ and $\alpha\in(\zeta,\pi/2)$, and let $ \{\bm{W}^t\}_{t \ge 0} \subseteq \St(d,k)$ be the sequence produced by Algorithm~\ref{alg:generic}. Suppose that Assumptions~\ref{as:eigengap} and \ref{as:delay} and the initialization condition \eqref{eq:init-zeta} hold, the stepsize 
	\[
		\eta \le \min \left\{
		\begin{aligned}
			\eta_0,\; 
			\frac{1}{4 \nu_{\rm avg} C_{\rm po} \sqrt{k}},\;
			\frac{1}{24 \nu_{\rm avg} (\tau+1)},\;
			\frac{1}{4\delta\cos^2(\alpha)(\tau+1)},\;
			\eta_{\rm comp}(\zeta,\alpha)
		\end{aligned}
		\right\}.
	\]
	Then, the iterates converge in Grassmannian geodesic distance as
	\begin{align}
		\dist_{\Gr}^2(\bm{W}^t,\bm{W}^*)
		 & \leq
		\left( 1 - \delta\cos^2(\alpha)\,\eta \right)^t
		\frac{\alpha^2(F(\bm{W}^0) - F^*)}{\delta\sin^2\alpha}.
		\label{eq:dist-linear}
	\end{align}
\end{theorem}

The constants in Theorem~\ref{thm} become more transparent for a concrete choice of radii.

\begin{remark}[Convergence under a concrete stepsize]\label{rem:zeta-gamma-choice}
	Take $\zeta={\pi}/{(3\sqrt2)}$ and $\alpha={\pi}/{3}$, then
	\(
		A_{\delta,\pi/3}={27\sqrt3\pi}/{\delta},
		B_{k,\pi/3}={17}/{8}+{2C_{\rm po}\sqrt{k}\pi}/{3},
		R_{\zeta,\alpha}={1}/{\nu},
	\)
	and the compatibility stepsize becomes
	\[
		\eta_{\rm comp}\left(\frac{\pi}{3\sqrt2},\frac{\pi}{3}\right)
		=
		\begin{cases}
			\displaystyle \frac{1}{4B_{k,\pi/3}\nu}, & \tau=0, \\[3mm]
			\displaystyle
			\min\left\{
			\frac{1}{4B_{k,\pi/3}\nu},\,
			\frac{1}{2\nu_{\rm avg}\tau\sqrt{144+27\sqrt3\pi\nu/\delta}}
			\right\}, & \tau\ge1.
		\end{cases}
	\]
	Under the stepsize conditions in Theorem~\ref{thm} together with this compatibility bound, the assumption $\dist_{\Gr}(\bm W^0,\bm W^*)\le\pi/(3\sqrt2)$ yields the concrete distance estimate
	\[
		\dist_{\Gr}^2(\bm W^t,\bm W^*)
		\le
		\left(1-\frac{\delta\eta}{4}\right)^t
		\frac{4\pi^2}{27\delta}
		\left(F(\bm W^0)-F^*\right).
	\]
\end{remark}

Solving the distance estimate \eqref{eq:dist-linear} with the largest admissible stepsize in Theorem~\ref{thm} gives the following iteration complexity.

\begin{corollary}[Broad-basin iteration complexity]\label{cor:broad-complexity}
	Under Theorem~\ref{thm}, choose the largest admissible stepsize allowed by the theorem, including the compatibility bound $\eta_{\rm comp}(\zeta,\alpha)$. Up to fixed constants depending on $\alpha$, $\zeta$, $k$, $\nu$, and $\nu_{\rm avg}$, the number of iterations $t$ sufficient for $\dist_{\Gr}^2(\bm W^t,\bm W^*)\le\varepsilon$ scales as
	\[
		\mathcal O\!\left(\delta^{-1}\log(1/(\delta\varepsilon))\right)
		\text{ when }\tau=0
		\quad\text{	and }\quad
		\mathcal O\!\left(\tau\delta^{-3/2}\log(1/(\delta\varepsilon))\right)\text{ when }\tau\ge1.
	\]
\end{corollary}

\subsubsection{Sharper Local Linear Convergence}\label{sec:sharper-local}

We next state sharper local guarantees under a stronger objective-gap initialization.

\begin{theorem}[Sharper local linear convergence]\label{thm:strong-init}
	Let $\{\bm W^t\}_{t\ge0}$ be the sequence produced by Algorithm~\ref{alg:generic}.
	Suppose that Assumptions~\ref{as:eigengap} and \ref{as:delay} hold.
	Assume the initialization condition
	\(
		F(\bm W^0)-F^* \le 4\delta/9.
	\)
	If the stepsize satisfies
	\[
		\eta \le \eta_{\rm loc}
		\coloneqq \min \left\{ \eta_0,\; \frac{1}{24\nu_{\rm avg}(\tau+1)},\; \frac{1}{\delta(\tau+1)} \right\},
	\]
	then all iterates remain in $\mathcal V_{\pi/3}$ and, for every $t\ge0$,
	\begin{align}
		\dist_{\Gr}^2(\bm W^t,\bm W^*)
		\le
		\left(1-\frac{\delta\eta}{4}\right)^t
		\frac{\pi^2(F(\bm W^0)-F^*)}{4\delta}.
		\label{eq:strong-init-dist}
	\end{align}
\end{theorem}

Solving the distance estimate \eqref{eq:strong-init-dist} with $\eta=\eta_{\rm loc}$ gives the following asymptotic iteration complexity in the local high-accuracy regime.

\begin{corollary}[Asymptotic iteration complexity]\label{cor:local-complexity}
	Under Theorem~\ref{thm:strong-init},
	choose
	\(
	\eta=\eta_{\rm loc}.
	\)
	Up to fixed constants depending on $k$, $\nu$, and $\nu_{\rm avg}$, the number of iterations $t$ sufficient for $\dist_{\Gr}^2(\bm W^t,\bm W^*)\le\varepsilon$ scales as
	\[
		\mathcal O\!\left( (\tau+1)\delta^{-1} \log(1/\varepsilon) \right),
	\]
\end{corollary}

Corollary~\ref{cor:local-complexity} describes the sharper asymptotic high-accuracy regime of the same {\algname} method.
After the iterates enter the stronger objective-gap region and the stepsize is set to the local bound, the extra broad-basin $\delta^{-1/2}$ penalty disappears: the delay dependence is linear and the eigengap dependence is $\mathcal O(\delta^{-1})$.

\paragraph{Shift-invariance of the results.}
It is worth noting that the convergence rates stated in Theorems~\ref{thm} and \ref{thm:strong-init} and the iteration complexities stated in Corollaries~\ref{cor:broad-complexity} and \ref{cor:local-complexity} are \emph{shift-invariant}: the constants depend on the component matrices only through the spectral spreads $\nu$ and $\nu_{\rm avg}$, rather than through absolute spectral levels such as $\max_i\|\bm A_i\|_2$.
This invariance removes a nuisance that is irrelevant to the eigenspace problem.
Indeed, componentwise shifts $\bm A_i\mapsto \bm A_i+c_i\bm I$ leave the target eigenspace and the eigengap unchanged, and they also leave the {\algname} trajectory unchanged because
\(
	(\bm I-\bm W\bm W^\top)(\bm A_i+c_i\bm I)\bm W
	=
	(\bm I-\bm W\bm W^\top)\bm A_i\bm W
	\text{ for all } c_i\in\mathbb R.
\)
Thus, the spread-based bound above is insensitive to this irrelevant parametrization, while a complexity bound depending on absolute spectral norms could be made arbitrarily worse by such shifts.

\section{Numerical Experiments}

We evaluate {\algname} on two finite-sum PCA tasks constructed from real datasets.
The first tests sample efficiency in the serial setting, and the second tests wall-clock convergence under distributed heterogeneous workers.
We report the objective gap $F(\bm W^t)-F(\bm W^*)$ or the geodesic distance $\dist_{\Gr}(\bm W^t,\bm W^*)$; shaded bands show one standard deviation over three independent runs.

\subsection{Serial \texorpdfstring{$k$}{k}-PCA}

\begin{figure}[t]
	\centering
	\includegraphics[width=1.0\textwidth]{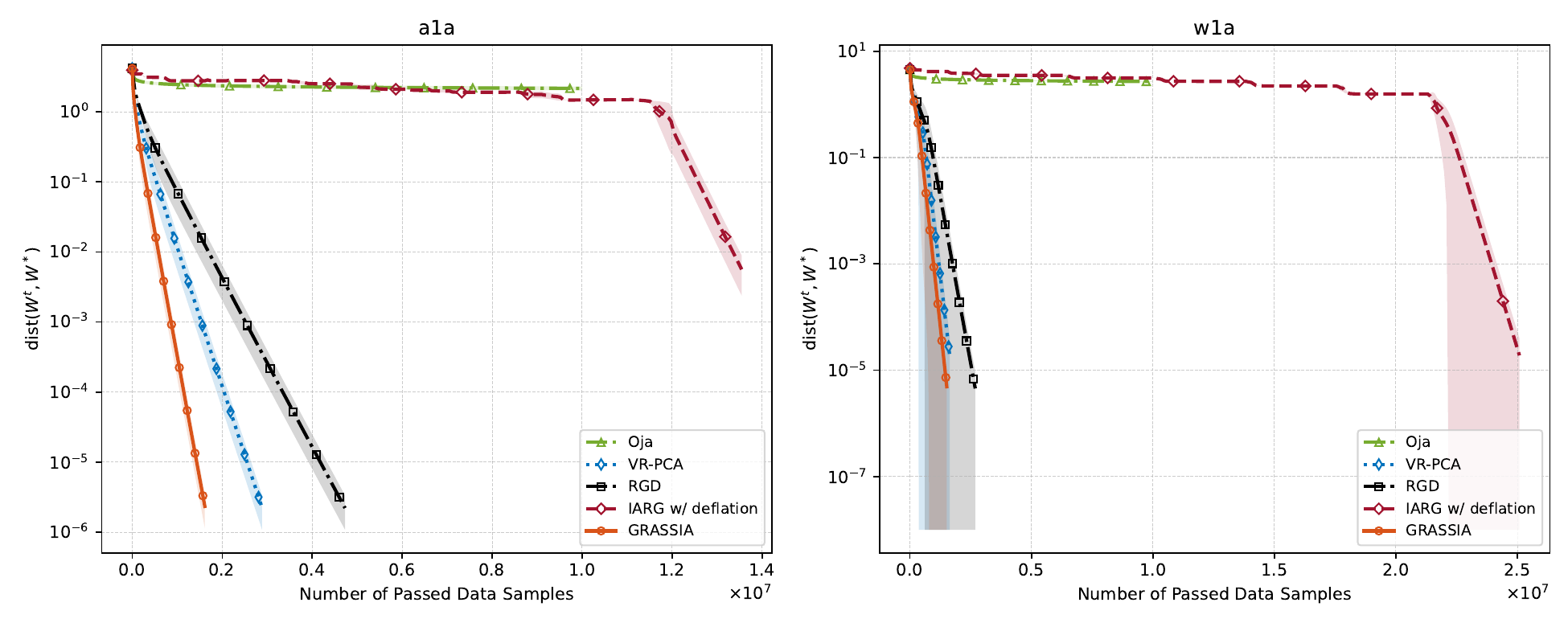}
	\caption{Incremental $k$-PCA on real datasets with $k=10$. Curves report the mean geodesic distance over three random initializations, and shaded bands report one standard deviation, with the lower band clipped only for log-scale visualization. The horizontal axis counts passed data samples. IARG with deflation solves each one-vector deflation stage to tolerance $10^{-6}$.}
	\label{fig:incremental-pca-real}
\end{figure}

We use the LIBSVM datasets \texttt{a1a} and \texttt{w1a} \citep{chang2011libsvm}.
For samples $\{\bm x_j\}_{j=1}^m\subset\mathbb R^d$, each component is the rank-one covariance contribution
\(
	\bm A_j=\bm x_j\bm x_j^\top
	\text{ so that }
	\bm A=\frac{1}{m}\sum_{j=1}^m \bm A_j.
\)
The instances have $(m,d)=(1605,123)$ for \texttt{a1a} and $(2477,300)$ for \texttt{w1a}; in both cases $k=10$, and $\bm W^*$ is computed from the full data matrix.
We compare GRASSIA with Block Oja, VR-PCA, RGD, and IARG with deflation.
The latter applies the single-vector IARG scheme sequentially: after one principal component is computed to tolerance $10^{-6}$, the covariance matrix is deflated and the next component is computed on the residual problem.
The incremental implementation uses mini-batches of $50$ samples, so the horizontal axis counts $50$ passed samples per update; one RGD iteration counts as one full data pass.
All methods use constant stepsizes selected by grid search. In the order Block Oja, VR-PCA, RGD, IARG with deflation, and GRASSIA, the stepsizes are
$(3\cdot10^{-2},\, 3\cdot10^{-2},\, 3\cdot10^{-1},\, 10^{-2},\, 3\cdot10^{-2})$ for \texttt{a1a}, and
$(3\cdot10^{-2},\, 5\cdot10^{-2},\, 8\cdot10^{-1},\, 10^{-2},\, 3\cdot10^{-2})$ for \texttt{w1a}.
Figure~\ref{fig:incremental-pca-real} shows that GRASSIA attains the smallest geodesic distance among the compared methods over the reported sample range on both datasets.
The deflation baseline and Oja remain at noticeably larger errors in the plotted range, while VR-PCA and RGD improve more slowly than GRASSIA.
This supports the roles of both aggregation and direct subspace optimization: aggregation reduces the noise of single-sample updates without requiring full-gradient iterations, and the direct $k$-dimensional update avoids the inefficiency of sequential deflation.

\begin{figure}[t]
	\centering
	\includegraphics[width=1.0\textwidth]{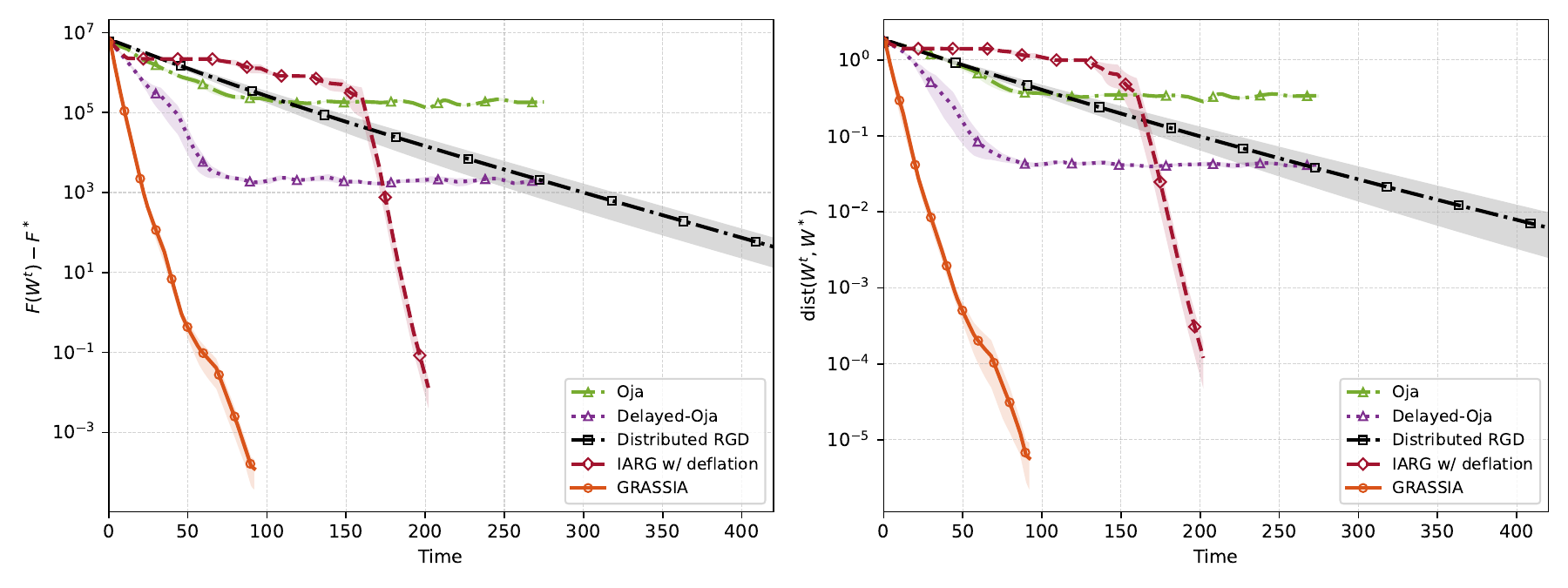}
	\caption{Asynchronous distributed PCA on CIFAR-10 with $d=3072$, $n=20$, and $k=3$. Curves include Oja, Delayed-Oja, Distributed RGD, IARG with deflation, and GRASSIA; all methods use the same data split and local batch size. Curves report the mean current value over three independent runs, and shaded bands report one half standard deviation. Left: objective gap versus time. Right: geodesic distance versus time.}
	\label{fig:async-pca-cifar}
\end{figure}

\subsection{Distributed \texorpdfstring{$k$}{k}-PCA with Heterogeneous Workers}

We next use $5000$ centered CIFAR-10 images \citep{krizhevsky2009learning}, represented in $\mathbb R^{3072}$ and evenly split across $n=20$ workers.
For the local data matrix $\bm X_i\in\mathbb R^{3072\times250}$ on worker $i$, set
\(
	\bm A_i=\frac{1}{250}\bm X_i \bm X_i^\top
	\text{ and }
	\bm A=\frac{1}{20}\sum_{i=1}^{20}\bm A_i.
\)
We set $k=3$.
To model heterogeneous computation, worker times $s_i$ are sampled uniformly from $\{1,2,3,4,5\}$; worker $i$ returns whenever the discrete clock is a multiple of $s_i$.
We compare with Oja, synchronous distributed RGD, and a non-aggregated stale-gradient ablation Delayed-Oja, whose update is
\(
	\bm W^{t+1}
	=
	\Polar \left( \bm W^t-\eta \grad F_{j_t}(\bm W^{t-\tau_{j_t}(t)}) \right).
\)
The RGD baseline waits for all workers each round and therefore takes $\max_i s_i$ clock ticks per full-gradient step.
All methods use grid-searched constant stepsizes fixed as
$(\eta_{\texttt{GRASSIA}},\eta_{\texttt{IARG}},\eta_{\texttt{Delayed-Oja}},\eta_{\texttt{RGD}},\eta_{\texttt{Oja}})=(10^{-7},2\cdot10^{-8},10^{-8},10^{-7},10^{-8})$.
Figure~\ref{fig:async-pca-cifar} shows that GRASSIA attains lower objective gaps and geodesic distances than the compared methods at the same wall-clock times over the reported range.
Distributed RGD also exhibits stable descent, but its synchronized communication rounds make it substantially slower under heterogeneous worker speeds.
Delayed-Oja and Oja stagnate at significantly larger error levels.
These observations indicate that stale-gradient aggregation is essential in the asynchronous setting: it preserves the low communication latency of incremental updates while using the cached local information to form a much more accurate search direction.

\section{Conclusion}

This paper develops {\algname}, a stale-aggregation method for finite-sum leading-eigenspace computation that keeps updates cheap by refreshing only arriving components.
The analysis proves shift-invariant two-phase linear convergence by coupling tight Grassmannian error bounds with delayed descent and basin invariance.
Experiments on serial and distributed PCA show consistent gains in sample efficiency and wall-clock convergence over baselines.
More broadly, the analysis suggests a framework for asynchronous Riemannian optimization that tailors stale-update control and basin invariance to the objective geometry, enabling problem-dependent guarantees across a wider class of manifold optimization problems.

\appendix

\section{Proofs of Results in Section~\ref{sec:geometry}}\label{app:geometry-proofs}

We first prove the two identities stated in Lemmas~\ref{lem:grad-energy} and \ref{lem:value-gap}, and then derive the angle-dependent gradient-dominance inequality.

\subsection{Proof of Lemma~\ref{lem:grad-energy}}
\begin{proof}
	Using the eigen-decomposition $\bm{A} = \bm{U}\bm{\Lambda} \bm{U}^\top$ and $\bm{Z} = \bm{U}^\top \bm{W}$, we write
	\begin{align}
		\frac{1}{4}\|\grad F(\bm{W})\|_F^2
		& =\|(\bm{I}-\bm{W}\bm{W}^\top)\bm{A}\bm{W}\|_F^2
		\nonumber
		\\
		& =\tr\left(\bm{W}^\top \bm{A}(\bm{I}-\bm{W}\bm{W}^\top)^2\bm{A} \bm{W}\right)
		\nonumber
		\\
		& =\tr(\bm{A}^2\bm{W}\bm{W}^\top) - \tr(\bm{A} \bm{W}\bm{W}^\top \bm{A} \bm{W}\bm{W}^\top)
		\nonumber
		\\
		& =\tr(\bm{\Lambda}^2 \bm{P})-\tr(\bm{\Lambda} \bm{P} \bm{\Lambda} \bm{P}).
		\nonumber
	\end{align}
	For the commutator expression, expand
	\begin{align}
		\|\comm{\bm{\Lambda}}{\bm{P}}\|_F^2
		= \| \bm{\Lambda} \bm{P}-\bm{P}\bm{\Lambda} \|_F^2
		=\tr(\bm{P}\bm{\Lambda}^2 \bm{P})-2 \tr(\bm{\Lambda} \bm{P}\bm{\Lambda} \bm{P})+\tr(\bm{\Lambda} \bm{P}^2\bm{\Lambda}).
		\nonumber
	\end{align}
	Because $\bm{P}^2=\bm{P}$, cyclicity gives $\tr(\bm{P}\bm{\Lambda}^2 \bm{P})=\tr(\bm{\Lambda}^2 \bm{P})$
	and $\tr(\bm{\Lambda} \bm{P}^2\bm{\Lambda})=\tr(\bm{\Lambda}^2 \bm{P})$, then we have
	$\|\comm{\bm{\Lambda}}{\bm{P}}\|_F^2=2\left(\tr(\bm{\Lambda}^2 \bm{P})-\tr(\bm{\Lambda} \bm{P} \bm{\Lambda} \bm{P})\right)$.
	Therefore,
	\begin{align}
		\frac{1}{4} \|\grad F(\bm{W})\|_F^2
		 & = \frac{1}{2} \|\comm{\bm{\Lambda}}{\bm{P}}\|_F^2.
		\label{eq:GE-comm}
	\end{align}
	Finally, block $\comm{\bm{\Lambda}}{\bm{P}}$ as
	\begin{align}
		\comm{\bm{\Lambda}}{\bm{P}}=
		\begin{bmatrix}
			\comm{\bm{\Lambda}_k}{\bm{P}_{11}}                             & \bm{\Lambda}_k \bm{P}_{12}-\bm{P}_{12}\bar{\bm{\Lambda}}_{d-k} \\
			\bar{\bm{\Lambda}}_{d-k} \bm{P}_{21}-\bm{P}_{21}\bm{\Lambda}_k & \comm{\bar{\bm{\Lambda}}_{d-k}}{\bm{P}_{22}}
		\end{bmatrix},
		\nonumber
	\end{align}
	sum the squared Frobenius norms of the four blocks, and apply $\bm{P}_{12}=\bm{P}_{21}^\top$ ($\bm{P}$ is symmetric) and \eqref{eq:GE-comm}, we have
	\begin{align}
		\frac{1}{4} \left\|\grad F(\bm{W})\right\|_F^2
		=\frac12 \left\|\comm{\bm{\Lambda}_k}{\bm{P}_{11}}\right\|_F^2
		+\frac12 \left\|\comm{\bar{\bm{\Lambda}}_{d-k}}{\bm{P}_{22}}\right\|_F^2
		+ \left\|\bm{\Lambda}_k\bm{P}_{12}\!-\!\bm{P}_{12}\bar{\bm{\Lambda}}_{d-k}\right\|_F^2.
		\nonumber
	\end{align}
	This completes the proof.
\end{proof}

\subsection{Proof of Lemma~\ref{lem:value-gap}}
\begin{proof}
	Using the eigen-decomposition $\bm{A} = \bm{U}\bm{\Lambda} \bm{U}^\top$ and $\bm{Z} = \bm{U}^\top \bm{W}$, we write
	\begin{align}
		F(\bm{W})
		= -\tr(\bm{W}^\top \bm{A} \bm{W})
		= -\tr(\bm{Z}^\top \bm{\Lambda} \bm{Z})
		= -\tr(\bm{Z}_1^\top\bm{\Lambda}_k \bm{Z}_1)-\tr(\bm{Z}_2^\top\bar{\bm{\Lambda}}_{d-k} \bm{Z}_2).
		\nonumber
	\end{align}
	Inserting the CS decomposition for both cases $d-k\ge k$ and $d-k<k$, the first block gives
	\begin{align}
		\tr(\bm{Z}_1^\top\bm{\Lambda}_k \bm{Z}_1)
		=\tr\left(
		\bm{R} \bm{C} \bm{Q}_1^\top \bm{\Lambda}_k \bm{Q}_1 \bm{C} \bm{R}^\top
		\right)
		=\tr\left(\bm{C} \bm{Q}_1^\top \bm{\Lambda}_k \bm{Q}_1 \bm{C}\right)
		=\tr(\bm{M}_k \bm{C}^2).
		\nonumber
	\end{align}
	For the second block, if $d-k\ge k$, set the row-blocked sine factor
	\(
		\bm S_{\mathrm{row}}
		\coloneqq
		\begin{bmatrix}
			\bm S \\
			\bm 0_{(d-2k)\times k}
		\end{bmatrix},
	\)
	so that $\bm Z_2=\bm Q_2\bm S_{\mathrm{row}}\bm R^\top$. Then,
	\begin{align}
		\tr(\bm{Z}_2^\top\bar{\bm{\Lambda}}_{d-k} \bm{Z}_2)
		& = \tr\left(
		\bm R\bm S_{\mathrm{row}}^\top
		\bm Q_2^\top\bar{\bm{\Lambda}}_{d-k}\bm Q_2
		\bm S_{\mathrm{row}}\bm R^\top
		\right)
		\nonumber
		\\
		& = \tr\left(
		\begin{bmatrix}
			\bm I_k & \bm 0_{k\times(d-2k)}
		\end{bmatrix}
		\bm Q_2^\top\bar{\bm{\Lambda}}_{d-k}\bm Q_2
		\begin{bmatrix}
			\bm I_k \\
			\bm 0_{(d-2k)\times k}
		\end{bmatrix}
		\bm S^2
		\right)
		\nonumber
		\\
		& = \tr(\overline{\bm M}_k\bm S^2).
		\nonumber
	\end{align}
	If $d-k<k$, set the column-blocked sine factor
	\(
	\bm S_{\mathrm{col}}\coloneqq
	\begin{bmatrix}
		\bm{0}_{(d-k)\times(2k-d)},\; \bm S_0
	\end{bmatrix},
	\)
	so that $\bm Z_2=\bm Q_2\bm S_{\mathrm{col}}\bm R^\top$ and
	\(
		\bm S_{\mathrm{col}}\bm S_{\mathrm{col}}^{\top}
		= \bm S_0^2.
	\)
	Then,
	\begin{align}
		\tr(\bm{Z}_2^\top\bar{\bm{\Lambda}}_{d-k} \bm{Z}_2)
		& = \tr\left(
		\bm R\bm S_{\mathrm{col}}^{\top}\bm Q_2^\top
		\bar{\bm{\Lambda}}_{d-k}\bm Q_2
		\bm S_{\mathrm{col}}\bm R^\top
		\right)
		\nonumber\\
		&= \tr\left(\bm Q_2^\top\bar{\bm{\Lambda}}_{d-k}\bm Q_2\bm S_0^2\right)
		\nonumber
		\\
		& = \tr\left(
		\begin{bmatrix}
				\lambda_{k+1}\bm I_{2k-d} & \bm 0                                       \\
				\bm 0                     & \bm Q_2^\top\bar{\bm{\Lambda}}_{d-k}\bm Q_2
			\end{bmatrix}
		\begin{bmatrix}
				\bm{0}_{2k-d} & \bm{0}    \\
				\bm{0}        & \bm S_0^2
			\end{bmatrix}
		\right)
		\nonumber
		\\
		& = \tr(\overline{\bm M}_k\bm S^2).
		\nonumber
	\end{align}
	Therefore,
	\begin{align}
		F(\bm{W})
		 & = - \tr(\bm{M}_k \bm{C}^2) - \tr(\overline{\bm{M}}_k \bm{S}^2).
		\label{eq:FW}
	\end{align}
	Besides, since $\bm{C}^2+\bm{S}^2=\bm{I}_k$ in both cases, we get
	\begin{align}
		F^*
		& = -\sum_{i=1}^k\lambda_i
		=-\tr(\bm{\Lambda}_k \bm{Q}_1 \bm{Q}_1^\top)
		=-\tr(\bm{M}_k)
		= -\tr(\bm{M}_k \bm{C}^2) - \tr(\bm{M}_k \bm{S}^2).
		\label{eq:F}
	\end{align}
	Combining \eqref{eq:FW} and \eqref{eq:F}, we have
	\(
	F(\bm{W})-F^* =\tr\left(\bm{M}_k \bm{S}^2\right) - \tr\left(\overline{\bm{M}}_k \bm{S}^2\right),
	\)
	which is the desired value-gap identity.
\end{proof}

\subsection{Proof of Proposition~\ref{lem:localPL}}
\begin{proof}
	By Lemma \ref{lem:grad-energy}, we have
	\begin{align}
		\frac{1}{4}\|\grad F(\bm{W})\|_F^2 \ge \|\bm{\Lambda}_k \bm{P}_{12}-\bm{P}_{12}\bar{\bm{\Lambda}}_{d-k}\|_F^2.
		\label{eq:offdiag-lb}
	\end{align}
	Let $\operatorname{vec}(\cdot)$ denote the column-wise vectorization operator. Then,
	\begin{align*}
		\operatorname{vec}(\bm{\Lambda}_k \bm{P}_{12})
		&=( \bm{I}_{d-k} \otimes \bm{\Lambda}_k ) \operatorname{vec}(\bm{P}_{12}),
		\\
		\operatorname{vec}(\bm{P}_{12}\bar{\bm{\Lambda}}_{d-k})
		&=\left( \bar{\bm{\Lambda}}_{d-k} \otimes \bm{I}_k \right) \operatorname{vec}(\bm{P}_{12}).
	\end{align*}
	Hence, we have
	\begin{align}
		\operatorname{vec}(\bm{\Lambda}_k \bm{P}_{12}-\bm{P}_{12}\bar{\bm{\Lambda}}_{d-k}) = \bm{\Omega} \operatorname{vec}(\bm{P}_{12}),
		\label{eq:vec-L}
	\end{align}
	where $\bm{\Omega} \coloneqq \bm{I}_{d-k} \otimes \bm{\Lambda}_k - \bar{\bm{\Lambda}}_{d-k} \otimes \bm{I}_k \in \mathbb{R}^{k(d-k)\times k(d-k)}$.
	Let $\bm e_i$ denote the standard basis vector (of appropriate dimension) whose $i$-th entry is 1 and all other entries are 0. For $i\in\{1,\dots,k\}$ and $j\in\{k+1,\dots,d\}$, define $\bm{\xi}_{(j,i)} \coloneqq \bm{e}_{j-k}\otimes \bm{e}_i \in \mathbb{R}^{k(d-k)}$ with $\bm{e}_{j-k} \in \mathbb{R}^{d-k}$ and $\bm{e}_i \in \mathbb{R}^k$. 
	Then, we have
	\begin{align}
		& (\bm{I}_{d-k} \otimes \bm{\Lambda}_k) \bm{\xi}_{(j,i)}
		=(\bm{I}_{d-k} \otimes \bm{\Lambda}_k) (\bm{e}_{j-k}\otimes \bm{e}_i)
		=\bm{e}_{j-k} \otimes (\lambda_i \bm{e}_i)
		= \lambda_i \bm{\xi}_{(j,i)},
		\nonumber
		\\
		& (\bar{\bm{\Lambda}}_{d-k} \otimes \bm{I}_k)\bm{\xi}_{(j,i)} =(\bar{\bm{\Lambda}}_{d-k} \otimes \bm{I}_k) (\bm{e}_{j-k}\otimes \bm{e}_i)
		=(\lambda_j \bm{e}_{j-k}) \otimes \bm{e}_i
		= \lambda_j \bm{\xi}_{(j,i)}.
		\nonumber
	\end{align}
	This implies that
	\(
		\bm{\Omega} \bm{\xi}_{(j,i)} = (\lambda_i-\lambda_j) \bm{\xi}_{(j,i)},
	\)
	i.e., each $\bm{\xi}_{(j,i)}$ is an eigenvector of $\bm{\Omega}$ with eigenvalue $\lambda_i-\lambda_j$.
	Note that the eigenvalues $\{\lambda_i-\lambda_j\}_{1 \le i \le k < j \le d}$ of $\bm{\Omega}$ satisfy $\lambda_i-\lambda_j \geq \lambda_k - \lambda_{k+1} = \delta$, then we have
	\(
	\bm{\Omega}^2 \succeq \delta \bm{\Omega}.
	\)
	This, together with \eqref{eq:vec-L}, implies that
	\begin{align}
		\left\| \bm{\Lambda}_k \bm{P}_{12}-\bm{P}_{12}\bar{\bm{\Lambda}}_{d-k} \right\|_F^2
		& = \left\| \bm{\Omega} \operatorname{vec}(\bm{P}_{12}) \right\|_2^2
		\nonumber
		\\
		& = \left\langle \operatorname{vec}(\bm{P}_{12}), \bm{\Omega}^2 \operatorname{vec}(\bm{P}_{12}) \right\rangle
		\nonumber
		\\
		& \ge \delta \langle \operatorname{vec}(\bm{P}_{12}), \bm{\Omega} \operatorname{vec}(\bm{P}_{12})\rangle
		\nonumber
		\\
		& = \delta \left( \left\langle \operatorname{vec}(\bm{P}_{12}), (\bm{I}_{d-k} \otimes \bm{\Lambda}_k) \operatorname{vec}(\bm{P}_{12}) \right\rangle
		- \left\langle \operatorname{vec}(\bm{P}_{12}), (\bar{\bm{\Lambda}}_{d-k} \otimes \bm{I}_k)\operatorname{vec}(\bm{P}_{12}) \right\rangle \right)
		\nonumber
		\\
		& = \delta \left( \tr(\bm{\Lambda}_k \bm{P}_{12} \bm{P}_{12}^\top)
		- \tr(\bar{\bm{\Lambda}}_{d-k} \bm{P}_{12}^\top \bm{P}_{12}) \right)
		\label{eq:diff-norm}
	\end{align}
	We next express the trace difference in \eqref{eq:diff-norm} using the case-wise CS decomposition in Lemma~\ref{lem:value-gap}.
	If $d-k\ge k$, then, with the row-blocked sine factor
	\(
	\bm S_{\mathrm{row}}
	=
	\begin{bmatrix}
		\bm S \\
		\bm 0_{(d-2k)\times k}
	\end{bmatrix},
	\)
	we have
	\[
		\bm{P}_{12} = \bm{Q}_1\bm{C}\bm S_{\mathrm{row}}^\top\bm{Q}_2^\top,
		\qquad
		\bm S_{\mathrm{row}}^\top\bm S_{\mathrm{row}} = \bm S^2.
	\]
	If $d-k<k$, then, with the column-blocked sine factor
	\(
		\bm S_{\mathrm{col}}=
		\begin{bmatrix}
			\bm{0}_{(d-k)\times(2k-d)},\; \bm{S}_0
		\end{bmatrix},
	\)
	we have
	\[
		\bm{P}_{12} = \bm{Q}_1\bm{C}\bm S_{\mathrm{col}}^{\top}\bm{Q}_2^\top,
		\qquad
		\bm S_{\mathrm{col}}^{\top}\bm S_{\mathrm{col}}=\bm S^2.
	\]
	Using the definitions of $\overline{\bm M}_k$ in both cases, one can verify that
	\[
		\tr(\bm{\Lambda}_k \bm{P}_{12} \bm{P}_{12}^\top) = \tr\left(\bm{M}_k \bm{C} \bm{S}^2 \bm{C}\right),
		\qquad
		\tr(\bar{\bm{\Lambda}}_{d-k}\bm{P}_{12}^\top \bm{P}_{12}) =
		\tr\left(\overline{\bm{M}}_k \bm{S} \bm{C}^2 \bm{S}\right).
	\]
	Thus, it follows that
	\begin{align*}
		\tr(\bm{\Lambda}_k \bm{P}_{12}\bm{P}_{12}^\top) - \tr(\bar{\bm{\Lambda}}_{d-k} \bm{P}_{12}^\top \bm{P}_{12})
		& = \tr\left((\bm{M}_k-\overline{\bm{M}}_k)  \bm{C}^2 \bm{S}^2 \right)
		\nonumber
		\\
		& \ge (\cos^2 \theta_k) \tr\left((\bm{M}_k-\overline{\bm{M}}_k)  \bm{S}^2\right)
		\nonumber
		\\
		& = (\cos^2 \theta_k) \left(F(\bm{W})-F^*\right),
	\end{align*}
	where the first equality is because $\bm{C}$ and $\bm{S}$ are both diagonal and $\bm{S} \bm{C}^2 \bm{S}=\bm{C} \bm{S}^2 \bm{C} = \bm{C}^2 \bm{S}^2$, the inequality holds due to $\bm{M}_k-\overline{\bm{M}}_k \succeq \bm 0$ and $\bm{C}^2 \bm{S}^2 \succeq (\cos^2 \theta_k) \bm{S}^2$ (since $\cos \theta_\ell \geq \cos \theta_k$ for all $\ell \in [k]$), and the last equality is the value-gap identity from Lemma \ref{lem:value-gap}.
	Inserting this into \eqref{eq:diff-norm} and combining with \eqref{eq:offdiag-lb} yield \eqref{eq:localPL} immediately.

	It remains to verify the tightness statement. 
	In the case $k=1$ and $d=2$, take $\bm A=\diag(\lambda_1,\lambda_2)$ with $\delta=\lambda_1-\lambda_2>0$, so that the leading eigenspace is spanned by $\bm e_1=[1,0]^\top$. For any unit vector $\bm w=(w_1,w_2)^\top$, let $\theta\in[0,\pi/2]$ be the principal angle between $\operatorname{span}(\bm w)$ and $\operatorname{span}(\bm e_1)$, thus $\cos\theta=|w_1|$ and $\sin\theta=|w_2|$. A direct computation gives
	\(
		F(\bm w)-F^*
		=
		\delta\sin^2\theta
		\text{ and }
		\|\grad F(\bm w)\|_F^2
		=
		4\delta^2\sin^2(\theta)\cos^2(\theta).
	\)
	Therefore,
	\[
		\frac12 \|\grad F(\bm w)\|_F^2 = 2\delta\cos^2(\theta)(F(\bm w)-F^*).
	\]
	Since $\delta>0$ and every $\theta\in[0,\pi/2]$ can be realized as the principal angle, any pointwise coefficient depending only on the eigengap and the principal angle that is strictly larger than $2\delta\cos^2\theta$ would violate the inequality on this two-dimensional instance.
\end{proof}

\section{Proofs of Results in Section~\ref{sec:one-step}}

This section proves the one-step estimates stated in Section~\ref{sec:one-step}.
The proof has two parts.
We first collect elementary stability bounds for the Riemannian gradients and for the aggregated direction; these bounds depend on the spectral spreads $\nu_i$ and are invariant under shifts of the component matrices.
We then use these bounds to show that the polar step is well defined, to expand its Grassmann logarithm up to second order, and finally to prove the descent inequality in Proposition~\ref{prop:one-step-descent}.

\begin{lemma}[Gradient Lipschitzness]\label{lem:smooth}
	It holds for all $\bm W,\bm W'\in\St(d,k)$ that
	\begin{align}
		\|\grad F_i(\bm W)-\grad F_i(\bm W')\|_F
		\le
		6\nu_i\|\bm W-\bm W'\|_F.
		\label{eq:smooth-Fi-2-F}
	\end{align}
\end{lemma}

\begin{proof}
	Let $\lambda_{d,i}\coloneqq\lambda_d(\bm A_i)$ and define $\widetilde{\bm A}_i\coloneqq\bm A_i-\lambda_{d,i}\bm I$. Since $(\bm I-\bm W\bm W^\top)\bm W=\bm 0$, we have the shift-invariance identity
	\[
		\grad F_i(\bm W)
		=-2(\bm I-\bm W\bm W^\top)\bm A_i\bm W
		=-2(\bm I-\bm W\bm W^\top)\widetilde{\bm A}_i\bm W.
	\]
	Moreover, $\|\widetilde{\bm A}_i\|_2=\lambda_1(\bm A_i)-\lambda_d(\bm A_i)=\nu_i$.
	Define projection matrix $\Pi(\bm W)\coloneqq\bm I-\bm W\bm W^\top$. Then,
	\begin{align*}
		\grad F_i(\bm W)-\grad F_i(\bm W')
		& = -2\Pi(\bm W)\widetilde{\bm A}_i\bm W + 2\Pi(\bm W')\widetilde{\bm A}_i\bm W' 
		\\
		& = -2\Pi(\bm W)\widetilde{\bm A}_i(\bm W-\bm W')
		-2(\Pi(\bm W)-\Pi(\bm W'))\widetilde{\bm A}_i\bm W'.
		\nonumber
	\end{align*}
	Note that $\|\Pi(\bm W)\|_2\le 1$ and
	\[
		\|\Pi(\bm W)-\Pi(\bm W')\|_F
		= \|(\bm W'-\bm W)\bm W'^\top+\bm W(\bm W'-\bm W)^\top \|_F
		\le 2\|\bm W-\bm W'\|_F
	\]
	due to $\|\bm W\|_2=\|\bm W'\|_2=1$.
	Then,
	\begin{align*}
		\|\grad F_i(\bm W)-\grad F_i(\bm W')\|_F
		&\le 2\|\Pi(\bm W)\|_2\|\widetilde{\bm A}_i\|_2\|\bm W-\bm W'\|_F
		\\
		&\quad + 2\|\Pi(\bm W)-\Pi(\bm W')\|_F
		\|\widetilde{\bm A}_i\|_2\|\bm W'\|_2
		\\
		&\le 2\nu_i\|\bm W-\bm W'\|_F + 2\cdot (2\|\bm W-\bm W'\|_F)\cdot \nu_i 
		\\
		&= 6\nu_i\|\bm W-\bm W'\|_F,
	\end{align*}
	which proves \eqref{eq:smooth-Fi-2-F}.
\end{proof}

\begin{lemma}[Aggregated gradient norm bound]\label{lem:D_G-explicit}
	For the aggregate Riemannian gradient \eqref{eq:approx-rg-gen}, one has the uniform bounds:
	\begin{align}
		\|\bm G^t\|_2 \le \nu_{\rm avg},
		\qquad
		\|\bm G^t\|_F \le \sqrt{k} \nu_{\rm avg}.
		\label{eq:G-uniform-bounds}
	\end{align}
\end{lemma}

\begin{proof}
	We first prove the corresponding bound for each component. Fix $i\in[n]$ and $\bm W\in\St(d,k)$. For any scalar $\alpha\in\mathbb R$,
	\(
	\grad F_i(\bm W) =-2(\bm I-\bm W\bm W^\top)(\bm A_i-\alpha\bm I)\bm W.
	\)
	Since $\|\bm I-\bm W\bm W^\top\|_2=1$ and $\|\bm W\|_2=1$, we have
	\begin{align}
		\|\grad F_i(\bm W)\|_2
		\le 2\|\bm I-\bm W\bm W^\top\|_2\|\bm A_i-\alpha\bm I\|_2\|\bm W\|_2
		=2\|\bm A_i-\alpha\bm I\|_2.
		\label{eq:grad-2norm}
	\end{align}
	We now minimize over $\alpha$. Because $\bm A_i$ is symmetric,
	\(
		\|\bm A_i-\alpha\bm I\|_2=\max\{|\lambda_1(\bm A_i)-\alpha|,\;|\lambda_d(\bm A_i)-\alpha|\}.
	\)
	The 1-dimensional problem
	\(
		\min_{\alpha\in\mathbb{R}} \max\{|\lambda_1(\bm A_i)-\alpha|,\;|\lambda_d(\bm A_i)-\alpha|\}
	\)
	attains its minimum at $\alpha_i^*=(\lambda_1(\bm A_i)+\lambda_d(\bm A_i))/2$. Thus,
	\[
		\min_{\alpha\in\mathbb R}\|\bm A_i-\alpha\bm I\|_2
		= \frac12 \left(\lambda_1(\bm A_i)-\lambda_d(\bm A_i)\right)
		= \frac{\nu_i}{2}.
	\]
	Substituting this into \eqref{eq:grad-2norm} gives
	\(
	\|\grad F_i(\bm W)\|_2 \le \nu_i
	\text{ and }
	\|\grad F_i(\bm W)\|_F \le \sqrt{k} \nu_i.
	\)
	Further using the aggregate definition \eqref{eq:approx-rg-gen}, we get
	\begin{align}
		\|\bm G^t\|_2
		 & \le \frac{1}{n} \sum_{i=1}^n \|\grad F_i(\bm W^{t-\tau_i(t)})\|_2
		\le \frac{1}{n} \sum_{i=1}^{n} \nu_i
		= \nu_{\rm avg}.
		\nonumber
	\end{align}
	Then, since $\bm G^t \in \mathbb R^{d\times k}$, we have $\|\bm G^t\|_F\le \sqrt{k} \|\bm G^t\|_2 \le \sqrt{k} \nu_{\rm avg}$.
\end{proof}

\begin{lemma}[Injectivity of the polar step]\label{lem:polar-inj-general}
	Let $\bm W^t\in\St(d,k)$ and let $\bm G^t\in\mathbb{R}^{d\times k}$ be defined as \eqref{eq:approx-rg-gen}.
	Assume that $0 < \eta < 1 / \nu_{\rm avg}$.
	Then, $\bm W^{t+1} \coloneqq \Polar(\bm W^t-\eta \bm G^t)$ and $\Log_{\bm W^t}(\bm W^{t+1})$ are both uniquely defined.
\end{lemma}

\begin{proof}
	Let $\bm X^t\coloneqq \bm W^t-\eta \bm G^t$.
	Using $\| \bm W^t \|_2 = 1$ and Lemma~\ref{lem:D_G-explicit}, we have
	\begin{align}
		\| (\bm W^t)^\top \bm G^t\|_2 \le \| \bm W^t \|_2 \| \bm G^t \|_2 \le \| \bm G^t \|_2 \le \nu_{\rm avg}.
		\nonumber
	\end{align}
	Thus, it follows from the step size condition $\eta < 1 / \nu_{\rm avg}$ that
	\(
	\eta \| (\bm W^t)^\top \bm G^t\|_2 < 1,
	\)
	which further implies that $\bm I_k-\eta (\bm W^t)^\top \bm G^t$ is invertible.
	If $\bm X^t\bm v=\bm 0$ for some $\bm v\in\mathbb{R}^k$, then premultiplying by $(\bm W^t)^\top$
	gives
	\(
	(\bm I_k-\eta (\bm W^t)^\top \bm G^t)\bm v=\bm 0,
	\)
	hence $\bm v=\bm 0$. Therefore, $\bm X^t$ has full column rank and $(\bm X^t)^\top \bm X^t\succ \bm 0$. This implies that $\bm W^{t+1} = \Polar(\bm X^t) = \bm X^t((\bm X^t)^\top\bm X^t)^{-1/2}$ is uniquely defined.
	Moreover, it follows that
	\begin{align}
		(\bm W^t)^\top \bm W^{t+1}
		= \left( (\bm W^t)^\top \bm X^t \right) \left((\bm X^t)^\top \bm X^t\right)^{-1/2}
		= \left( \bm I_k-\eta (\bm W^t)^\top \bm G^t \right) \left( (\bm X^t)^\top \bm X^t \right)^{-1/2},
		\nonumber
	\end{align}
	where both factors are invertible. Thus, $(\bm W^t)^\top \bm W^{t+1}$ is invertible, and consequently
	\begin{align*}
		\sigma_{\min}((\bm W^t)^\top\bm W^{t+1})&>0,
		\\
		\cos\theta_i(\bm W^t,\bm W^{t+1})
		&=\sigma_i((\bm W^t)^\top\bm W^{t+1})>0,
		\qquad i\in[k].
	\end{align*}
	Equivalently,
	\(
	\theta_k(\bm W^t,\bm W^{t+1})<\pi/2.
	\)
	By the logarithm formula in Section~\ref{sec:notation}, this places $\bm W^{t+1}$ in the injectivity domain of $\Log_{\bm W^t}(\cdot)$, and hence $\Log_{\bm W^t}(\bm W^{t+1})$ is uniquely defined.
\end{proof}

\begin{lemma}[Upper bounds on geodesic distance]\label{lem:geodesic-step}
	Let $\bm{W}^t$ and $\bm{W}^{t+1}$ with $t \geq 0$ be any two consecutive iterates produced by Algorithm~\ref{alg:generic}. If the stepsize satisfy
	\(
	\eta \le 1 / (2\nu_{\rm avg}),
	\)
	then the geodesic step on $\Gr(d,k)$ for the polar update satisfies
	\begin{align}
		\dist_{\Gr}(\bm W^{t+1},\bm W^t)
		\le 2\eta \|\bm G^t\|_F
		\le
		2\eta\sqrt{k}\nu_{\rm avg}.
		\label{eq:dist-step-explicit}
	\end{align}
\end{lemma}

\begin{proof}
	Let $\bm X^t \coloneqq \bm W^t-\eta\bm G^t$.
	Since
	\(
		\eta \le 1/(2 \nu_{\rm avg}) < 1/\nu_{\rm avg},
	\)
	it follows from Lemma \ref{lem:polar-inj-general} that $\bm W^{t+1} = \Polar(\bm X^t)$ is uniquely defined.
	Then, we have
	\begin{align}
		&\dist_{\Gr}(\bm W^{t+1},\bm W^t)
		= \|\bm{\Theta}(\bm W^{t+1},\bm W^t)\|_F
		\leq \|\tan\bm{\Theta}(\bm W^{t+1},\bm W^t)\|_F
		= \|\tan\bm{\Theta}(\bm X^t,\bm W^t)\|_F.
		\label{eq:dist-tan}
	\end{align}
	where the inequality uses the fact that $\theta\le \tan\theta$ for $\theta \in [0,\pi/2)$ and the last equality is because $\mathrm{span}(\bm W^{t+1})=\mathrm{span}(\bm X^t)$.
	Let $\bm{\sigma}(\cdot)$ be the vector of singular values and $\bm W^t_\perp\in\mathbb R^{d\times(d-k)}$ be any orthonormal complement of $\bm W^t$ so that $[\bm W^t,\bm W^t_\perp]\in\mathbb O(d)$.
	The tangent principal-angle formula \citep{zhu2013angles} gives
	\[
		\tan\bm{\Theta}(\bm X^t,\bm W^t) = \bm{\sigma} ((\bm W^t_\perp)^\top \bm X^t ((\bm W^t)^\top \bm X^t)^{-1}).
	\]
	Compute  
	\(
	(\bm W^t_\perp)^\top \bm X^t = -\eta(\bm W^t_\perp)^\top \bm G^t
	\text{ and }
	(\bm W^t)^\top \bm X^t = \bm I_k-\eta(\bm W^t)^\top \bm G^t,
	\)
	then we have
	\begin{align}
		\|\tan\bm{\Theta}(\bm W^{t+1},\bm W^t)\|_F
		&= \left\| \bm{\sigma} \left(\eta(\bm W^t_\perp)^\top \bm G^t (\bm I_k-\eta(\bm W^t)^\top \bm G^t)^{-1} \right)\right\|_2
		\nonumber
		\\
		&= \left\| (\eta(\bm W^t_\perp)^\top \bm G^t (\bm I_k-\eta(\bm W^t)^\top \bm G^t)^{-1} \right\|_F
		\nonumber
		\\
		&\le \eta\|\bm G^t\|_F \|(\bm I_k-\eta(\bm W^t)^\top \bm G^t)^{-1}\|_2,
		\nonumber
	\end{align}
	where the inequality is bacause $\|(\bm W^t_\perp)^\top \bm G^t\|_F = \|(\bm I-\bm W^t(\bm W^t)^\top)\bm G^t\|_F\le \|\bm G^t\|_F$.
	By $\| \bm W^t \|_2 = 1$ and Lemma~\ref{lem:D_G-explicit},
	\(
	\|\eta (\bm W^t)^\top\bm G^t\|_2
	\le \eta \|\bm W^t\|_2\|\bm G^t\|_2
	\le \eta \nu_{\rm avg} 
	\le 1/2.
	\)
	This yields the Neumann-series bound 
	\(
		\|(\bm I_k-\eta(\bm W^t)^\top \bm G^t)^{-1}\|_2
		\le (1-\eta\|(\bm W^t)^\top \bm G^t\|_2)^{-1}.
	\)
	Therefore, it holds that
	\[
	\|\tan\bm{\Theta}(\bm W^{t+1},\bm W^t)\|_F
	\le \eta \|\bm G^t\|_F \cdot \frac{1}{1-\eta\|(\bm W^t)^\top \bm G^t\|_2}
	\le 2 \eta \|\bm G^t\|_F
	\leq 2\eta\sqrt{k}\nu_{\rm avg},
	\]
	where the last inequality uses Lemma~\ref{lem:D_G-explicit}.
	Substituting this into \eqref{eq:dist-tan} gives the desired bound \eqref{eq:dist-step-explicit}.
\end{proof}

\subsection{Proof of Lemma~\ref{lem:log-decomp}}
The proof below rewrites the polar update in affine coordinates around the current subspace and then compares this coordinate representation with the Grassmann logarithm.
The resulting remainder is cubic and is dominated by the stated second-order bound under the stepsize condition.

\begin{proof}
	Since stepsize condition \eqref{eq:eta-cond} implies that
	\(
		\eta \le 1/(2 \nu_{\rm avg}) < 1/\nu_{\rm avg},
	\)
	it follows from Lemma \ref{lem:polar-inj-general} that $\bm W^{t+1} = \Polar(\bm W^t-\eta \bm G^t)$ and $\Log_{\bm W^t}(\bm W^{t+1})$ are both uniquely defined.
	Let $\bm{G}_{\perp}^t \coloneqq \mathscr{P}_{\mathcal{T}_{[\bm W^t]}\Gr(d,k)}(\bm G^t) = (\bm I-\bm W^t(\bm W^t)^\top)\bm G^t$, then we write $\Log_{\bm W^t}(\bm W^{t+1}) = - \eta\bm{G}_{\perp}^t + \bm E^t$. Then, it suffices to upper bound the Frobenius norm of the error $\bm E^t = \Log_{\bm W^t}(\bm W^{t+1}) + \eta\bm{G}_{\perp}^t$.
	We split
	\begin{align}
		\bm E^t
		=(\bm L^t+\eta\bm{G}_{\perp}^t)
		+(\Log_{\bm W^t}(\bm W^{t+1})-\bm L^t).
		\label{eq:error-split}
	\end{align}
	Define
	\[
		\bm{X}^t\coloneqq \bm W^t-\eta\bm G^t,
		\qquad
		\bm{H}^t\coloneqq (\bm W^t)^\top \bm G^t.
	\]
	Then, we have
	\(
	\bm W^{t+1}=\bm{X}^t((\bm{X}^t)^\top\bm{X}^t)^{-1/2},
	\)
	\(
	(\bm W^t)^\top\bm W^{t+1}
	=(\bm I-\eta \bm{H}^t)((\bm{X}^t)^\top\bm{X}^t)^{-1/2},
	\)
	and
	\begin{align}
		\bm L^t
		\coloneqq&\; \bm W^{t+1}\big((\bm W^t)^\top\bm W^{t+1}\big)^{-1}-\bm W^t 
		\nonumber\\
		=&\; \bm{X}^t(\bm I-\eta\bm{H}^t)^{-1}-\bm W^t 
		\nonumber\\
		=&\; (\bm W^t-\eta\bm G^t)(\bm I-\eta\bm{H}^t)^{-1}
		-\bm W^t (\bm I-\eta\bm{H}^t) (\bm I-\eta\bm{H}^t)^{-1}
		\nonumber\\
		=&\;-\eta(\bm G^t-\bm W^t\bm{H}^t)(\bm I-\eta\bm{H}^t)^{-1}
		\nonumber\\
		=&\;-\eta\bm{G}_{\perp}^t (\bm I-\eta\bm{H}^t)^{-1}.
		\nonumber
	\end{align}
	Therefore, the first two terms in \eqref{eq:error-split} can be written as
	\begin{align}
		&\bm E_1^t
		\coloneqq \bm L^t+\eta\bm{G}_{\perp}^t 
		= -\eta\bm{G}_{\perp}^t(\bm I-\eta\bm{H}^t)^{-1} + \eta\bm{G}_{\perp}^t 
		= -\eta \bm{G}_{\perp}^t (\bm I-\eta\bm{H}^t)^{-1} \eta \bm{H}^t, 
		\nonumber
	\end{align}
	which implies that
	\begin{align}
		\|\bm E_1^t\|_F
		&\le \eta^2\|\bm G_{\perp}^t\|_F\|\bm H^t\|_2\|(\bm I-\eta\bm H^t)^{-1}\|_2
		\nonumber\\
		&\le \eta^2\|\bm G_{\perp}^t\|_F\|\bm H^t\|_2 \cdot \frac{1}{1-\eta\|\bm H^t\|_2} 
		\nonumber\\
		&\le \eta^2\|\bm G^t\|_F^2 \cdot \frac{1}{1-\eta\|\bm H^t\|_2} 
		\nonumber\\
		&\le 2\eta^2\|\bm G^t\|_F^2, 
		\label{eq:E1-bound}
	\end{align}
	where the second inequality used the Neumann-series bound for $\eta\|\bm H^t\|_2\le 1/2$ (due to the stepsize condition \eqref{eq:eta-cond}), the thrid inequality is because $\|\bm G_{\perp}^t\|_F=\|(\bm I-\bm W^t(\bm W^t)^\top)\bm G^t\|_F\le \|\bm G^t\|_F$, and the last inequality follows from $\|\bm H^t\|_2=\|(\bm W^t)^\top\bm G^t\|_2\le  \|\bm G^t\|_2\le\nu_{\rm avg}$.
	Let the thin SVD of $\bm L^t$ be $\bm L^t=\bm R\bm\Sigma\bm V^\top$ with $\bm R\in\St(d,k)$, $\bm V\in\mathbb O(k)$, and $\bm\Sigma$ diagonal, then the Grassmann logarithm is given by
	\(
		\Log_{\bm W^t}(\bm W^{t+1})=\bm R\arctan(\bm\Sigma)\bm V^\top,
	\)
	Thus, the last two terms in \eqref{eq:error-split} can be written as
	\begin{align}
		\bm E_2^t
		 & \coloneqq \Log_{\bm W^t}(\bm W^{t+1})-\bm L^t
		= \bm R\big(\arctan(\bm\Sigma)-\bm\Sigma\big)\bm V^\top, \nonumber
	\end{align}
	Therefore, we have
	\begin{align}
		\|\bm E_2^t\|_F
		&= \| \arctan(\bm\Sigma)-\bm\Sigma \|_F
		\nonumber
		\\
		&\le \frac{1}{3} \|\bm\Sigma\|_F^3
		= \frac{1}{3} \|\bm R \bm\Sigma \bm V^\top\|_F^3
		= \frac{1}{3} \|\bm L^t\|_F^3
		= \frac{1}{3}\|\bm L^t\|_F^2\cdot \|\bm L^t\|_F
		\nonumber
		\\
		&\le \frac{1}{3}\cdot (2\eta\|\bm G^t\|_F)^2 \cdot \sqrt{k}
		\nonumber
		\\
		&\le \frac{4\sqrt{k}}{3}\eta^2\|\bm G^t\|_F^2,
		\label{eq:E2-bound}
	\end{align}
	where the first inequality is because
	\[
		0\le x-\arctan(x)=\int_{0}^{x}\frac{u^2}{1+u^2}{\rm d}u \le \int_{0}^{x}u^2 {\rm d} u=\frac{x^3}{3}
	\]
	for $x\ge 0$, the second inequality follows from
	\[
		\|\bm L^t\|_F
		\le \eta\|\bm G^t\|_F\|(\bm I-\eta\bm{H}^t)^{-1}\|_2
		\le 2\eta\|\bm G^t\|_F
		\le \sqrt{k},
	\]
	and the last inequality uses Lemma~\ref{lem:D_G-explicit} and the stepsize condition \eqref{eq:eta-cond}.
	Finally, substituting \eqref{eq:E1-bound} and \eqref{eq:E2-bound} into \eqref{eq:error-split} yields
	\[
		\|\bm E^t\|_F \le \|\bm E_2^t\|_F+\|\bm E_1^t\|_F
		\le C_{\rm po} \eta^2 \|\bm G^t\|_F^2,
	\]
	which proves \eqref{eq:log-decomp}.
\end{proof}

\subsection{Proof of Proposition~\ref{prop:one-step-descent}}
The proof combines the logarithmic expansion from Lemma~\ref{lem:log-decomp} with a second-order Taylor bound for $F$ on $\Gr(d,k)$.

\begin{proof}
	By Lemma \ref{lem:log-decomp}, the stepsize condition $\eta \le 1/(2\nu_{\rm avg})$ guarantees that $\bm W^{t+1} = \Polar(\bm W^t-\eta \bm G^t)$ and $\Log_{\bm W^t}(\bm W^{t+1})$ are both well defined.
	Then, it follows from the second-order Taylor expansion of $F$ on $\Gr(d,k)$ \citep{boumal2023introduction} and Lemma \ref{lem:geodesic-step} that
	\begin{align}
		F(\bm{W}^{t+1}) - F(\bm{W}^{t})
		\leq & \left\langle \grad F(\bm{W}^{t}), \Log_{\bm{W}^t}(\bm{W}^{t+1}) \right\rangle
		+ \nu \dist_{\Gr}^2(\bm{W}^{t+1}, \bm{W}^{t})
		\nonumber
		\\
		\leq & \left\langle \grad F(\bm{W}^{t}), -\eta \mathscr{P}_{\mathcal{T}_{[\bm W^t]} \Gr(d,k)}(\bm G^t) + \bm{E}^t \right\rangle
		+ 4 \nu \eta^2 \|\bm G^t\|_F^2.
		\label{eq:fval-diff}
	\end{align}
	Since
	\(
		\bm{G}^t
		= \mathscr{P}_{\mathcal{T}_{[\bm{W}^t]} \Gr(d,k)} (\bm{G}^t)
		+ \mathscr{N}_{\mathcal{T}_{[\bm{W}^t]} \Gr(d,k)} (\bm{G}^t)
	\)
	and 
	$\grad F(\bm{W}^{t}) \in \mathcal{T}_{[\bm{W}^t]} \Gr(d,k)$,
	we have
	\[
		\langle \grad F(\bm{W}^{t}), \bm{G}^t \rangle
		= \langle \grad F(\bm{W}^{t}), \mathscr{P}_{\mathcal{T}_{[\bm{W}^t]} \Gr(d,k)} (\bm{G}^t) \rangle.
	\]
	Applying this to inequality \eqref{eq:fval-diff} and using Lemma \ref{lem:log-decomp}, we get
	\begin{align}
		F(\bm{W}^{t+1}) - F(\bm{W}^{t})
		\leq & \left\langle \grad F(\bm{W}^{t}), -\eta \bm G^t \right\rangle
		+ \left\langle \grad F(\bm{W}^{t}), \bm{E}^t \right\rangle
		+ 4 \nu \eta^2 \|\bm G^t\|_F^2
		\nonumber
		\\
		\le  & - \eta \left\langle \grad F(\bm{W}^{t}), \bm{G}^t \right\rangle
		+ C_{\rm po} \sqrt{k} \nu \eta^2\|\bm G^t\|_F^2
		+ 4 \nu \eta^2 \|\bm G^t\|_F^2.
		\label{eq:fval-diff1}
	\end{align}
	Note that polar decomposition is the Euclidean nearest-point projection onto $\St(d,k)$, i.e.,
	\(
		\bm{W}^{s+1} = \argmin_{\bm{W} \in \St(d,k)} \| \bm{W} - (\bm{W}^s - \eta \bm{G}^s) \|_F
	\)
	for all $s \in \mathbb{N}$.
	This implies that
	\[
	\left\| \bm{W}^{s+1} - \bm{W}^s \right \|_F - \left \| \eta \bm{G}^s \right\|_F
	\leq \left\| \bm{W}^{s+1} - (\bm{W}^s - \eta \bm{G}^s) \right\|_F
	\le \eta \| \bm{G}^s \|_F.
	\]
	Thus, we have
	\(
	\left\| \bm{W}^{s+1} - \bm{W}^s \right\|_F
		\leq 2 \eta \| \bm{G}^s \|_F.
	\)
	This, together with the telescopic sum identity, implies that for any $i \in [n]$ and $t\ge0$ with $\tau_i(t) \ge 1$,
	\begin{align}
		\left\| \bm{W}^t - \bm{W}^{t-\tau_i(t)} \right\|_F
		\le \sum_{s = t-\tau_i(t)}^{t-1} \left\| \bm{W}^{s+1} - \bm{W}^s  \right\|_F
		\leq 2 \eta \sum_{s = [t-\tau]_+}^{t-1} \| \bm{G}^s \|_F. \label{eq:iter-diff}
	\end{align}
	Note that for any $i \in [n]$ and $t\ge0$ satisfying $\tau_i(t) = 0$, the bound \eqref{eq:iter-diff} holds trivially, since in this case $\| \bm{W}^t - \bm{W}^{t-\tau_i(t)} \|_F^2 = 0$.
	Let 
	\[
		\Delta \bm G^t \coloneqq \bm G^t - \grad F(\bm W^t) = \frac{1}{n} \sum_{i=1}^{n} \left( \grad F_i( \bm{W}^{t-\tau_i(t)} ) - \grad F_i( \bm{W}^t ) \right).
	\]
	Then, it follows from Lemma \ref{lem:smooth} that
	\begin{align}
		&\|\Delta \bm G^t\|_F
		\le \frac{6}{n} \sum_{i=1}^{n} \nu_i \| \bm{W}^{t-\tau_i(t)} - \bm{W}^t \|_F
		\le 12 \nu_{\rm avg} \eta \sum_{s = [t-\tau]_+}^{t-1} \| \bm{G}^s \|_F.
		\nonumber
	\end{align}
	Squaring it and applying Cauchy-Schwarz inequality, we have
	\begin{align}
		\|\Delta \bm G^t\|_F^2
		 & \le (12 \nu_{\rm avg} \eta)^2 \left( \sum_{s = [t-\tau]_+}^{t-1} \| \bm{G}^s \|_F \right)^2
		\le 144 \nu_{\rm avg}^2 \eta^2 \tau \sum_{s = [t-\tau]_+}^{t-1} \| \bm{G}^s \|_F^2\label{eq:DeltaG-upper}
	\end{align}
	Thus, the inner product term in \eqref{eq:fval-diff1} can be lower bounded as
	\begin{align}
		\left\langle \grad F(\bm{W}^{t}), \bm{G}^t \right\rangle
		&= \frac{1}{2} \| \grad F (\bm{W}^{t}) \|_F^2
		+ \frac{1}{2} \| \bm{G}^t \|_F^2
		- \frac{1}{2} \| \Delta \bm G^t \|_F^2
		\nonumber
		\\
		&\geq \frac{1}{2} \| \grad F (\bm{W}^{t}) \|_F^2
		+ \frac{1}{2} \| \bm{G}^t \|_F^2
		- 72 \nu_{\rm avg}^2 \tau \eta^2 \sum_{s = [t-\tau]_+}^{t-1} \| \bm{G}^s \|_F^2.
		\nonumber
	\end{align}
	Besides, it holds for all $t\ge0$ that
	\begin{align}
		\|\bm{G}^t\|_F^2
		& = \| \Delta \bm G^t + \grad F(\bm W^t) \|_F^2
		\leq 288 \nu_{\rm avg}^2 \tau \eta^2 \sum_{s = [t-\tau]_+}^{t-1} \| \bm{G}^s \|_F^2
		+ 2 \| \grad F (\bm{W}^{t}) \|_F^2.
		\nonumber 
	\end{align} 
	Substituting the above two bounds into \eqref{eq:fval-diff1}, we have
	\begin{align}
		& F(\bm{W}^{t+1}) - F(\bm{W}^{t})
		\leq - \left( \frac{1}{2} \eta - 2 (C_{\rm po} \sqrt{k} + 4)\nu \eta^2 \right) \| \grad F (\bm{W}^{t}) \|_F^2
		- \frac{1}{2} \eta \| \bm{G}^t \|_F^2
		\nonumber\\
		&\quad\; + \left(
		72 \nu_{\rm avg}^2 \tau \eta^3
		+ 288 (C_{\rm po} \sqrt{k} + 4)
		\nu \nu_{\rm avg}^2 \tau \eta^4
		\right)
		\sum_{s=[t-\tau]_+}^{t-1} \| \bm{G}^s \|_F^2
		\nonumber\\
		&\leq - \frac{1}{4} \eta \| \grad F (\bm{W}^{t}) \|_F^2
		- \frac{1}{2} \eta \| \bm{G}^t \|_F^2
		+ 144\nu_{\rm avg}^2\tau\eta^3
		\sum_{s=[t-\tau]_+}^{t-1}\!\| \bm{G}^s \|_F^2,
		\nonumber
	\end{align}
	where the second inequality holds because $\eta \leq 1/(8\nu (C_{\rm po} \sqrt{k} + 4))$ implies that 
	\(
		{\eta}/{2} - 2 (C_{\rm po} \sqrt{k} + 4)\nu \eta^2 \geq {\eta}/{4}
		\text{ and }
		72 \nu_{\rm avg}^2 \tau \eta^3 + 288 (C_{\rm po} \sqrt{k} + 4) \nu \nu_{\rm avg}^2 \tau \eta^4 \leq 144 \nu_{\rm avg}^2 \tau \eta^3.
	\)
	This completes the proof.
\end{proof}

\section{Proofs of Results in Section~\ref{sec:basin-invariance}}

This section proves the basin-invariance result used in the convergence theorem.
The proof combines an objective-gap bound, a Grassmann cosine-law estimate, and the one-step decrease from Proposition~\ref{prop:one-step-descent} in a Lyapunov argument.

The following bounds can also be obtained by the smoothness and quadratic-growth results of \citet[Propositions~3 and 6]{alimisis2024geodesic} for points with $\theta_k(\bm W,\bm W^*)<\pi/2$.
We give a concise independent proof based on Lemma~\ref{lem:value-gap}.
Since the proof uses only the value-gap identity, it does not require the logarithm map to be uniquely defined and therefore also covers the boundary case $\theta_k(\bm W,\bm W^*)=\pi/2$.

\begin{lemma}[Objective-gap bounds]\label{prop:obj-bound-gr}
	Suppose that Assumption \ref{as:eigengap} holds. Then, for any $\bm W\in\St(d,k)$, it holds that
	\begin{align}
		({4\delta}/{\pi^2}) \dist_{\Gr}^2(\bm W,\bm W^*)
		\le
		F(\bm W)-F^*
		\le
		\nu \dist_{\Gr}^2(\bm W,\bm W^*).
		\label{eq:gap-geodesic}
	\end{align}
\end{lemma}

\begin{proof}
	From the definitions of $\bm M_k$ and $\overline{\bm M}_k$ in both cases $d-k \ge k$ and $d-k < k$, we have
	$\lambda_{\min}(\bm M_k) = \lambda_k$,
	$\lambda_{\max}(\bm M_k) = \lambda_1$,
	$\lambda_{\min}(\overline{\bm M}_k) \ge \lambda_d$,
	and
	$\lambda_{\max}(\overline{\bm M}_k) \le \lambda_{k+1}$.
	Hence,
	\[
		\bm M_k-\overline{\bm M}_k
		\succeq (\lambda_k-\lambda_{k+1})\bm I_k=\delta \bm I_k,
		\qquad
		\bm M_k-\overline{\bm M}_k
		\preceq (\lambda_1-\lambda_d)\bm I_k=\nu \bm I_k.
	\]
	Combining this with $\bm S^2\succeq \bm0$ and Lemma~\ref{lem:value-gap} yields
	\begin{align}
		\delta \tr(\bm S^2)
		\le \tr \big((\bm M_k-\overline{\bm M}_k)\bm S^2\big)
		\le \nu \tr(\bm S^2),
		\nonumber
	\end{align}
	which is equivalent to
	\begin{align}
		\delta\sum_{i=1}^k\sin^2\theta_i
		\le
		F(\bm W)-F^*
		\le
		\nu\sum_{i=1}^k\sin^2\theta_i.
		\label{eq:gap-chordal}
	\end{align}
	Using the fact that $(2/\pi)\theta \le \sin\theta \le \theta$ for $\theta\in[0,\pi/2]$, we have
	\[
		\frac{4}{\pi^2}\sum_{i=1}^k \theta_i^2
		\le \sum_{i=1}^k \sin^2\theta_i
		\le \sum_{i=1}^k \theta_i^2,
	\]
	which is equivalent to
	\[
		\frac{4}{\pi^2} \dist_{\Gr}^2(\bm W,\bm W^*)
		\le \|\sin\bm{\Theta}(\bm W,\bm W^*)\|_F^2
		\le \dist_{\Gr}^2(\bm W,\bm W^*).
	\]
	Combining this with \eqref{eq:gap-chordal} yields \eqref{eq:gap-geodesic}.
\end{proof}

Then, we introduce the following cosine-law inequality on $\Gr(d,k)$ (equipped with the canonical metric), which is the trigonometric distance bound of \citet[Lemma~1]{zhang2016riemannian} specialized to curvature lower bound zero.
This bound follows from Toponogov's triangle comparison theorem \citep{CheegerEbin1975Comparison}; the fact that $\Gr(d,k)$ has nonnegative sectional curvature is standard; see, e.g., \citet{bendokat2024grassmann} for a concise account.

\begin{lemma}[Cosine-law inequality on $\Gr(d,k)$]\label{lem:3-point}
	Let $\bm X,\bm Y,\bm Z \in \Gr(d,k)$.
	Suppose that both $\Log_{\bm X}(\bm Y)$ and $\Log_{\bm X}(\bm Z)$ are uniquely defined. Then,
	\[
	\dist_{\Gr}^2(\bm Z,\bm Y)
	\le \|\Log_{\bm X}(\bm Y)\|_F^2 + \|\Log_{\bm X}(\bm Z)\|_F^2
	-2\left\langle \Log_{\bm X}(\bm Z), \Log_{\bm X}(\bm Y) \right\rangle.
	\]
\end{lemma}

\subsection{Proof of Proposition~\ref{prop:max-angle}}
The proof is an induction on the iteration count.
Assuming that the previous iterates lie in $\mathcal V_\alpha$, we show that the Lyapunov function does not increase up to time $T$. The explicit bound $\eta\le\eta_{\rm comp}(\zeta,\alpha)$ then converts this into the desired distance bound $\dist_{\Gr}(\bm W^T,\bm W^*)\le\alpha$.

\begin{proof}
	We establish the Proposition via induction. The base case follows from
	\(
		\dist_{\Gr}(\bm W^0,\bm W^*) \le \zeta < \alpha.
	\)
	Therefore, $\theta_k(\bm W^0,\bm W^*)\le\alpha$, and equivalently $\bm W^0\in\mathcal V_\alpha$.
	As the induction hypothesis, it holds for some integer $T\ge1$ that
	\[
		\dist_{\Gr}(\bm W^t,\bm W^*)\le\alpha,
		\qquad t=0,\ldots,T-1.
	\]
	For any $t \in \{0,\dots,T-1\}$, since $\theta_k(\bm W^t,\bm W^*)\le\dist_{\Gr}(\bm W^t,\bm W^*)\le\alpha<\pi/2$ and $\theta_k(\bm W^t,\bm W^{t+1})<\pi/2$ (by Lemma \ref{lem:polar-inj-general}), both $\Log_{\bm W^t}(\bm W^*)$ and $\Log_{\bm W^t}(\bm W^{t+1})$ are uniquely defined.
	Then, it follows Lemma~\ref{lem:3-point}, identity $\| \Log_{\bm W^t}(\bm W^*) \|_F^2 = \dist_{\Gr}^2(\bm W^t,\bm W^*)$, and Lemma~\ref{lem:log-decomp} that
	\begin{align}
		&\dist_{\Gr}^2(\bm W^{t+1},\bm W^*)-\dist_{\Gr}^2(\bm W^t,\bm W^*) 
		\nonumber
		\\
		\leq\;& \| \Log_{\bm W^t}(\bm W^{t+1}) \|_F^2
		-2 \left\langle
		\Log_{\bm W^t}(\bm W^{t+1}),
		\Log_{\bm W^t}(\bm W^*) \right\rangle
		\nonumber
		\\
		\le\;& 2\eta^2\|\bm G^t\|_F^2
		+ 2\|\bm E^t\|_F^2 + 2\eta \left\langle \bm G^t, \Log_{\bm W^t}(\bm W^*) \right\rangle
		+ 2\|\bm E^t\|_F \|\Log_{\bm W^t}(\bm W^*) \|_F.
		\label{eq:dist-part-raw}
	\end{align}
	By Lemmas~\ref{lem:log-decomp} and \ref{lem:D_G-explicit}, and $\eta \le 1/(4 \nu_{\rm avg} C_{\rm po} \sqrt{k})$, we have
	\[
		\|\bm E^t\|_F
		\le C_{\rm po} \eta^2\|\bm G^t\|_F^2
		\le C_{\rm po} \sqrt{k} \nu_{\rm avg} \eta^2 \|\bm G^t\|_F
		\le \frac{\eta}{4}\|\bm G^t\|_F,
	\]
	which implies that
	\begin{align}
		& 2\|\bm E^t\|_F^2 \le 2\cdot\frac{\eta^2}{16}\|\bm G^t\|_F^2=\frac{\eta^2}{8}\|\bm G^t\|_F^2.
		\label{eq:E^2-bounds}
	\end{align}
	Moreover, by the induction hypothesis, $\theta_i(\bm W^t,\bm W^*)\le \theta_k(\bm W^t,\bm W^*)\le \alpha$ for all $i\in[k]$, and hence
	\(
		\|\Log_{\bm W^t}(\bm W^*)\|_F=\dist_{\Gr}(\bm W^t,\bm W^*)\le \sqrt{k}\alpha.
	\)
	Therefore,
	\begin{align}
		& 2\|\bm E^t\|_F\; \|\Log_{\bm W^t}(\bm W^*)\|_F
		\le 2 C_{\rm po} \eta^2\|\bm G^t\|_F^2 \cdot \sqrt{k}\alpha
		= 2C_{\rm po}\sqrt{k}\alpha\, \eta^2\|\bm G^t\|_F^2.
		\label{eq:r-bounds}
	\end{align}
	We treat the inner-product
	$\langle \bm G^t,\Log_{\bm W^t}(\bm W^*)\rangle$ by splitting
	$\bm G^t=\grad F(\bm W^t)+\Delta \bm G^t$:
	\begin{align}
		\langle \bm G^t,\Log_{\bm W^t}(\bm W^*)\rangle
		&= \left\langle \grad F(\bm W^t),\Log_{\bm W^t}(\bm W^*) \right\rangle
		+ \left\langle \Delta \bm G^t,\Log_{\bm W^t}(\bm W^*) \right\rangle \nonumber                   \\
		&\le - \frac{2\theta_k(\bm W^t,\bm W^*)}{\tan\theta_k(\bm W^t,\bm W^*)} \left( F(\bm W^t)-F^* \right)
		+ \left\langle \Delta \bm G^t,\Log_{\bm W^t}(\bm W^*) \right\rangle  \nonumber                   \\
		& \le - \frac{8\delta}{\pi^2} \frac{\theta_k(\bm W^t,\bm W^*)}{\tan\theta_k(\bm W^t,\bm W^*)} \dist_{\Gr}^2(\bm W^t,\bm W^*) + \|\Delta \bm G^t\|_F \dist_{\Gr}(\bm W^t,\bm W^*), \nonumber
	\end{align}
	where the first inequality is due to the weak-quasi-convexity of Problem \eqref{eq:eigen} for $\bm W^t$ satisfying $\theta_k(\bm W^t,\bm W^*)<\pi/2$ \citep[Proposition~3]{alimisis2024geodesic} and the second inequality uses Lemma~\ref{prop:obj-bound-gr} and identity $\|\Log_{\bm W^t}(\bm W^*)\|_F = \dist_{\Gr}(\bm W^t,\bm W^*)$.
	Since $\theta/\tan\theta$ is monotonically decreasing on $(0,\pi/2)$, then it follows from the induction hypothesis that
	\[
		\frac{\theta_k(\bm W^t,\bm W^*)}{\tan\theta_k(\bm W^t,\bm W^*)}
		\ge \frac{\alpha}{\tan \alpha}.
	\]
	Substituting this and using Young's inequality $2 x y \le {x^2}/{\beta} + \beta y^2$ with $\beta = {16\delta\alpha}/(\pi^2\tan\alpha)$, we have
	\begin{align}
		2\eta \langle \bm G^t,\Log_{\bm W^t}(\bm W^*)\rangle
		&\le - \frac{16\delta\alpha}{\pi^2\tan\alpha} \eta \dist_{\Gr}^2(\bm W^t\!,\!\bm W^*)
		\!+\! \frac{\pi^2\tan\alpha}{16\delta\alpha} \eta \|\Delta \bm G^t\|_F^2
		\!+\! \frac{16\delta\alpha}{\pi^2\tan\alpha} \eta \dist_{\Gr}^2(\bm W^t\!,\!\bm W^*)
		\nonumber\\
		&= \frac{\pi^2\tan\alpha}{16\delta\alpha} \eta \|\Delta \bm G^t\|_F^2
		\nonumber\\
		&\le A_{\delta,\alpha} \nu_{\rm avg}^2 \tau \eta^3 \sum_{s = [t-\tau]_+}^{t-1} \| \bm{G}^s \|_F^2,
		\nonumber
	\end{align}
	where the last inequality uses \eqref{eq:DeltaG-upper}.
	Substituting this bound, together with \eqref{eq:E^2-bounds} and \eqref{eq:r-bounds}, into \eqref{eq:dist-part-raw} gives
	\begin{align}
		& \dist_{\Gr}^2(\bm W^{t+1},\bm W^*)-\dist_{\Gr}^2(\bm W^t,\bm W^*)
		\le B_{k,\alpha}\eta^2\|\bm G^t\|_F^2
		+ A_{\delta,\alpha} \nu_{\rm avg}^2 \tau \eta^3 \sum_{s = [t-\tau]_+}^{t-1} \!\!\| \bm{G}^s \|_F^2.
		\nonumber
	\end{align}
	From Proposition~\ref{prop:one-step-descent}, dropping the negative term
	$-\frac14\eta\|\grad F(\bm W^t)\|_F^2$ yields
	\begin{align}
		F(\bm W^{t+1})-F(\bm W^t)
		\le
		- \frac12\eta\|\bm G^t\|_F^2
		+ 144 \nu_{\rm avg}^2 \tau \eta^3 \sum_{s=[t-\tau]_+}^{t-1} \| \bm{G}^s \|_F^2.
		\nonumber
	\end{align}
	Define the Lyapunov function for $t \in \mathbb{N}$:
	\[
	\mathcal L^t \coloneqq
	\gamma_{\alpha,\eta}\left(F(\bm W^t)-F^*\right)
	+ \dist_{\Gr}^2(\bm W^t,\bm W^*),
	\] 
	where 
	\[
		\gamma_{\alpha,\eta}
		\coloneqq
		\begin{cases}
			\displaystyle
			\max\left\{
			4B_{k,\alpha}\eta,\,
			\frac{4 A_{\delta,\alpha}\nu_{\rm avg}^2\tau^2\eta^2}
			{1-576\nu_{\rm avg}^2\tau^2\eta^2}
			\right\}, & \tau\ge1,   
			\\[4mm]
			4B_{k,\alpha}\eta, & \tau=0.
		\end{cases}
	\]
	Then, using $\gamma_{\alpha,\eta}\ge 4B_{k,\alpha}\eta$, it holds for $t=0,\ldots,T-1$ that 
	\begin{align}
		\mathcal L^{t+1}-\mathcal L^{t}
		&= \gamma_{\alpha,\eta}\left( F(\bm W^{t+1})-F(\bm W^t) \right)
		+ \left( \dist_{\Gr}^2(\bm W^{t+1},\bm W^*)-\dist_{\Gr}^2(\bm W^t,\bm W^*) \right) \nonumber
		\\
		&\le - \frac{\gamma_{\alpha,\eta}}{4} \eta \|\bm G^t\|_F^2
		+ \left(144\gamma_{\alpha,\eta}+A_{\delta,\alpha}\right)
		\nu_{\rm avg}^2\tau\eta^3
		\sum_{s=[t-\tau]_+}^{t-1}\|\bm G^s\|_F^2.
		\nonumber
	\end{align}
	Summing this inequality over $t=0,\ldots,T-1$ and applying the elementary counting estimate
	\(
	\sum_{t=0}^{T-1}\sum_{s=[t-\tau]_+}^{t-1}\|\bm G^s\|_F^2
	\le
	\tau\sum_{t=0}^{T-1}\|\bm G^t\|_F^2,
	\)
	we get
	\begin{align}
		\mathcal L^{T}-\mathcal L^{0}
		\le -\left( \frac{\gamma_{\alpha,\eta}}{4} \eta
		- \left(144\gamma_{\alpha,\eta}+A_{\delta,\alpha}\right) \nu_{\rm avg}^2 \tau^2 \eta^3 \right)
		\sum_{t=0}^{T-1}\|\bm G^t\|_F^2
		\le 0.
		\nonumber
	\end{align}
	The second inequality in the preceding display is immediate if $\tau=0$.  If $\tau\ge1$, the stepsize condition gives
	$1-576\nu_{\rm avg}^2\tau^2\eta^2>0$ and the definition of $\gamma_{\alpha,\eta}$ implies
	\(
		\gamma_{\alpha,\eta}
		\left({1}/{4}-144\nu_{\rm avg}^2\tau^2\eta^2\right)
		\ge
		A_{\delta,\alpha}\nu_{\rm avg}^2\tau^2\eta^2,
	\)
	which is equivalent to the coefficient in parentheses being nonnegative.
	Thus, we have
	\begin{align}
		& \dist_{\Gr}^2(\bm W^T,\bm W^*)
		\le \mathcal L^T
		\le \mathcal L^0
		\le (1+\nu\gamma_{\alpha,\eta}) \dist_{\Gr}^2(\bm W^0,\bm W^*),
		\nonumber
	\end{align}
	where the last inequality uses Lemma~\ref{prop:obj-bound-gr}. 
	Then, we record the consequence
	\begin{align}
		(1+\nu\gamma_{\alpha,\eta})\zeta^2\le\alpha^2.
		\label{eq:zeta-alpha-compat}
	\end{align}
	This is immediate if $\zeta=0$.  If $\zeta>0$, then, by the definition of
	$R_{\zeta,\alpha}$, \eqref{eq:zeta-alpha-compat} is equivalent to
	$\gamma_{\alpha,\eta}\le R_{\zeta,\alpha}$.
	When $\tau\ge1$, $\eta \le \eta_{\rm comp}(\zeta,\alpha) \le R_{\zeta,\alpha}/(4B_{k,\alpha})$ gives
	\(4B_{k,\alpha}\eta\le R_{\zeta,\alpha}\), and 
	\[
		\eta \le \eta_{\rm comp}(\zeta,\alpha) \le \frac{1}{\nu_{\rm avg}\tau} \sqrt{\frac{R_{\zeta,\alpha}}{4A_{\delta,\alpha}+576R_{\zeta,\alpha}}}
		\;\Longrightarrow\;
		\frac{4 A_{\delta,\alpha}\nu_{\rm avg}^2\tau^2\eta^2}{1-576\nu_{\rm avg}^2\tau^2\eta^2}
		\le R_{\zeta,\alpha}.
	\]
	This, according to the definition of $\gamma_{\alpha,\eta}$ for $\tau \ge 1$, implies that $\gamma_{\alpha,\eta}\le R_{\zeta,\alpha}$.
	When $\tau=0$, the bound $\eta\le\eta_{\rm comp}(\zeta,\alpha) = R_{\zeta,\alpha}/(4B_{k,\alpha})$ immediately gives
	\(
	\gamma_{\alpha,\eta}=4B_{k,\alpha}\eta\le R_{\zeta,\alpha},
	\)
	proving \eqref{eq:zeta-alpha-compat}.
	Finally, it follows from \eqref{eq:init-zeta} and \eqref{eq:zeta-alpha-compat} that
	\[
		\dist_{\Gr}(\bm W^T,\bm W^*)
		\le \sqrt{1+\nu\gamma_{\alpha,\eta}} \cdot \zeta
		\le \alpha.
	\]
	Thus, $\theta_k(\bm W^T,\bm W^*)\le \alpha$, i.e., $\bm W^T \in \mathcal{V}_{\alpha}$.
	The induction step is proved.
\end{proof}

\section{Proofs of Results in Section~\ref{sec:two-phase}}

This section turns the descent and basin-invariance results into the two convergence theorems.
The following lemma, due to \citet[Lemma~1]{aytekin2016analysis}, will be used to absorb the delayed tail term in the resulting objective-gap recursion.
\begin{lemma}\label{lem:recur}
	Let $\{V_t\}_{t\ge0}$ and $\{W_t\}_{t\ge0}$ be nonnegative real sequences satisfying
	\(
	V_{t+1} \leq a V_t - b W_t + c \sum_{s=(t-T_0)_+}^{t} W_s
	\)
	for all $t\in\mathbb{N}$, where $a \in (0,1)$, $b,c \geq 0$, and $T_0\in\mathbb{N}$.
	Suppose that
	\[
	\begin{aligned}
		\frac{c}{1-a} \frac{1-a^{T_0+1}}{a^{T_0}} \leq b.
	\end{aligned}
	\]
	Then, $V_t \leq a^t V_0$ for all $t\in\mathbb{N}$.
\end{lemma}

\subsection{Proof of Theorem~\ref{thm}}
\begin{proof}
	Set $\Delta_t\coloneqq F(\bm W^t)-F^*$.
	Since the initial point satisfies \eqref{eq:init-zeta} and the stepsize satisfies the bounds required in Proposition~\ref{prop:max-angle}, that proposition implies that $\bm W^t\in\mathcal V_\alpha$ for all $t\ge0$.
	Then, by Proposition \ref{lem:localPL}, it holds for $t\ge0$ that
	\begin{align}
		\| \grad F (\bm{W}^t) \|_F^2
		\geq 4 \delta \cos^2 (\theta_k) \Delta_t
		\geq 4\delta\cos^2(\alpha) \Delta_t.
		\nonumber
	\end{align}
	Combining this with Proposition \ref{prop:one-step-descent} yields the recursion:
	\begin{align}
		\Delta_{t+1}
		\leq\left( 1 - \delta\cos^2(\alpha)\,\eta \right) \Delta_t
		- \frac{1}{2} \eta \| \bm{G}^t \|_F^2
		+ 144 \nu_{\rm avg}^2 \tau \eta^3 \sum_{s=[t-\tau]_+}^{t-1} \| \bm{G}^s \|_F^2&,
		\; t \geq 0.
		\label{eq:recur-Delta}
	\end{align}
	To apply Lemma \ref{lem:recur} with $T_0=\tau$, we let $a = 1 - \delta\cos^2(\alpha)\,\eta$, $b = \frac{1}{2} \eta$, and $c = L^2 \tau \eta^3$ with $L \coloneqq 12 \nu_{\rm avg}$. If $\tau=0$, then $c=0$ and the condition in Lemma~\ref{lem:recur} holds trivially. It remains to consider $\tau\ge1$.
	The stepsize condition $\eta\le 1/(4\delta\cos^2(\alpha)(\tau+1))$ gives
	\(
		1-a=\delta\cos^2(\alpha)\,\eta \le 1/(4(\tau+1))
		\text{ and thus }
		a^\tau\ge 1-\tau(1-a)\ge 3/4.
	\)
	Therefore,
	\[
		\frac{c}{1-a}\frac{1-a^{\tau+1}}{a^\tau}
		= \frac{c}{a^\tau} \sum_{\ell=0}^{\tau}a^\ell
		\le \frac43 (\tau+1) c
		\le 2L^2\tau(\tau+1)\eta^3.
	\]
	The sufficient condition
	\(
		\eta \le 1/(24\nu_{\rm avg}(\tau+1))
	\)
	ensures \(2L^2\tau(\tau+1)\eta^3\le\eta/2=b\).
	Also, $\eta\le1/(4\delta\cos^2(\alpha)(\tau+1))$ implies $a\in(0,1)$.
	Therefore, applying Lemma \ref{lem:recur} to \eqref{eq:recur-Delta} gives
	\begin{align}
		\Delta_t
		\le
		\left( 1 - \delta\cos^2(\alpha)\,\eta \right)^t
		\Delta_0,
		\qquad t\ge0.
		\label{eq:fval-linear}
	\end{align}
	For the distance estimate, Proposition~\ref{prop:max-angle} gives $\theta_i(\bm W^t,\bm W^*)\le\alpha$ for all $i\in[k]$.
	Using \eqref{eq:gap-chordal} and the fact that $x \mapsto \sin x / x$ is monitonically non-decreasing, we have
	\[
		\Delta_t
		\ge
		\delta\sum_{i=1}^k\sin^2\theta_i(\bm W^t,\bm W^*)
		\ge
		\delta ({\sin^2(\alpha)}/{\alpha^2})
		\dist_{\Gr}^2(\bm W^t,\bm W^*).
	\]
	Combining this bound with \eqref{eq:fval-linear} yields \eqref{eq:dist-linear}.
\end{proof}

\subsection{Proof of Theorem~\ref{thm:strong-init}}

\begin{proof}
	We argue by induction that the iterates stay in $\mathcal V_{\pi/3}$; the proof uses an auxiliary objective-gap recursion to close the induction.
	By Lemma~\ref{prop:obj-bound-gr} and the initialization condition,
	\[
		\dist_{\Gr}^2(\bm W^0,\bm W^*)
		\le
		\frac{\pi^2}{4\delta} \Delta_0
		\le
		\frac{\pi^2}{9}.
	\]
	Hence, $\theta_k(\bm W^0,\bm W^*)\le\dist_{\Gr}(\bm W^0,\bm W^*)\le \pi/3$, and therefore $\bm W^0\in\mathcal V_{\pi/3}$.
	Assume that $\bm W^t\in\mathcal V_{\pi/3}$ for all $t=0,\ldots,T$.
	Inside $\mathcal V_{\pi/3}$ we have $\cos^2\theta_k(\bm W^t,\bm W^*)\ge1/4$, and hence $\|\grad F(\bm W^t)\|_F^2\ge\delta\Delta_t$ by Proposition~\ref{lem:localPL}.
	Then, Proposition~\ref{prop:one-step-descent} imply that
	\begin{align}
		\Delta_{t+1}
		& \le
		\left(1-\frac{\delta\eta}{4}\right)\Delta_t
		-\frac{\eta}{2}\|\bm G^t\|_F^2
		+144\nu_{\rm avg}^2\tau\eta^3
		\sum_{r=[t-\tau]_+}^{t-1}\|\bm G^r\|_F^2,
		\quad t=0,\ldots,T.
		\label{eq:strong-local-rec}
	\end{align}
	We apply the Lemma~\ref{lem:recur} with
	\(
		a=1-{\delta\eta}/{4},\,
		b={\eta}/{2},\,
		c=(12\nu_{\rm avg})^2\tau\eta^3,
		\text{ and }
		T_0=\tau.
	\)
	The condition in Lemma~\ref{lem:recur} follows from the same verification as in the proof of Theorem~\ref{thm}, with $\alpha=\pi/3$, since the stepsize bounds reduce to $\eta\le1/(\delta(\tau+1))$ and $\eta\le1/(24\nu_{\rm avg}(\tau+1))$.
	Hence Lemma~\ref{lem:recur} yields, on the induction horizon,
	\[
		\Delta_t
		\le
		\left(1-\frac{\delta\eta}{4}\right)^t\Delta_0,
		\qquad t=0,\ldots,T+1.
	\]
	In particular, 
	\(
		\Delta_{T+1}\le\Delta_0\le4\delta/9.
	\)
	Using Lemma~\ref{prop:obj-bound-gr} once more gives 
	\[
		\dist_{\Gr}^2(\bm W^{T+1},\bm W^*)
		\le
		\frac{\pi^2}{4\delta}\Delta_{T+1}
		\le
		\frac{\pi^2}{9},
	\]
	so $\bm W^{T+1}\in\mathcal V_{\pi/3}$.
	This closes the induction. Therefore, 
	\[
		\Delta_t
		\le
		\left(1-\frac{\delta\eta}{4}\right)^t\Delta_0,
		\quad t=0,1,2,\ldots.
	\]
	The distance estimate \eqref{eq:strong-init-dist} follows by combining it with Lemma~\ref{prop:obj-bound-gr}.
\end{proof}

\bibliographystyle{arxiv-numeric}
\bibliography{ref}

@inproceedings{zhang2016riemannian,
	title={Riemannian {SVRG}: Fast stochastic optimization on {Riemannian} manifolds},
	author={Zhang, Hongyi and Reddi, Sashank J. and Sra, Suvrit},
	booktitle={Adv. Neural Inf. Process. Syst. 29},
	year={2016}
}

@article{guo2024fedpower,
	title   = {{FedPower}: Privacy-preserving Distributed Eigenspace Estimation},
	author  = {Guo, Xiao and Li, Xiang and Chang, Xiangyu and Wang, Shusen and Zhang, Zhihua},
	journal = {Mach. Learn.},
	volume  = {113},
	pages   = {8427--8458},
	year    = {2024}
}

@article{ye2021deepca,
	title={{DeEPCA}: Decentralized Exact {PCA} with Linear Convergence Rate.},
	author={Ye, Haishan and Zhang, Tong},
	journal={J. Mach. Learn. Res.},
	volume={22},
	number={238},
	pages={1--27},
	year={2021}
}

@inproceedings{huang2020communication,
	title={Communication-efficient distributed {PCA} by {R}iemannian optimization},
	author={Huang, Long-Kai and Pan, Sinno},
	booktitle={Proc. 37th Int. Conf. Mach. Learn.},
	pages={4465--4474},
	year={2020},
	organization={PMLR}
}

@inproceedings{alimisis2021distributed,
	title={Distributed Principal Component Analysis with Limited Communication},
	author={Alimisis, Foivos and Davies, Peter and Vandereycken, Bart and Alistarh, Dan},
	booktitle={Adv. Neural Inf. Process. Syst. 34},
	pages={2823--2834},
	year={2021}
}

@inproceedings{chen2021decentralized,
	title={Decentralized {R}iemannian Gradient Descent on the {S}tiefel Manifold},
	author={Chen, Shixiang and Garcia, Alfredo and Hong, Mingyi and Shahrampour, Shahin},
	booktitle={Proc. 38th Int. Conf. Mach. Learn.},
	pages={1594--1605},
	year={2021},
	organization={PMLR}
}

@inproceedings{li2021communication,
	title={Communication-efficient distributed {SVD} via local power iterations},
	author={Li, Xiang and Wang, Shusen and Chen, Kun and Zhang, Zhihua},
	booktitle={Proc. 38th Int. Conf. Mach. Learn.},
	pages={6504--6514},
	year={2021},
	organization={PMLR}
}

@article{gang2021distributed,
	title={Distributed principal subspace analysis for partitioned big data: Algorithms, analysis, and implementation},
	author={Gang, Arpita and Xiang, Bingqing and Bajwa, Waheed U},
	journal={IEEE Trans. Signal Inf. Process. Netw.},
	volume={7},
	pages={699--715},
	year={2021},
	publisher={IEEE}
}

@article{andrade2023distributed,
	title={Distributed {B}anach-{P}icard Iteration: Application to Distributed Parameter Estimation and {PCA}},
	author={Andrade, Francisco L. and Figueiredo, M{\'a}rio A. T. and Xavier, Jo{\~a}o},
	journal={IEEE Trans. Signal Process.},
	volume={71},
	pages={17--30},
	year={2023},
	publisher={IEEE}
}

@inproceedings{grammenos2020federated,
	title={Federated principal component analysis},
	author={Grammenos, Andreas and Mendoza Smith, Rodrigo and Crowcroft, Jon and Mascolo, Cecilia},
	booktitle={Adv. Neural Inf. Process. Syst. 33},
	pages={6453--6464},
	year={2020}
}

@inproceedings{garber2016faster,
	title={Faster eigenvector computation via shift-and-invert preconditioning},
	author={Garber, Dan and Hazan, Elad and Jin, Chi and Kakade, Sham M. and Musco, Cameron and Netrapalli, Praneeth and Sidford, Aaron},
	booktitle={Proc. 33rd Int. Conf. Mach. Learn.},
	pages={2626--2634},
	year={2016},
	organization={PMLR}
}

@inproceedings{xu2018accelerated,
	title={Accelerated stochastic power iteration},
	author={Xu, Peng and He, Bryan and De Sa, Christopher and Mitliagkas, Ioannis and Re, Chris},
	booktitle={Proc. 21st Int. Conf. Artif. Intell. Stat.},
	pages={58--67},
	year={2018},
	organization={PMLR}
}

@inproceedings{kim2020stochastic,
	title={Stochastic variance-reduced algorithms for {PCA} with arbitrary mini-batch sizes},
	author={Kim, Cheolmin and Klabjan, Diego},
	booktitle={Proc. Int. Conf. Artif. Intell. Stat.},
	pages={4302--4312},
	year={2020},
	organization={PMLR}
}

@article{zhu2013angles,
	title={Angles between subspaces and their tangents},
	author={Zhu, Peizhen and Knyazev, Andrew V.},
	journal={J. Numer. Math.},
	volume={21},
	number={4},
	pages={325--340},
	year={2013}
}

@inproceedings{tang2019matrixkrasulina,
	title     = {Exponentially Convergent Stochastic $k$-{PCA} without Variance Reduction},
	author    = {Tang, Cheng},
	booktitle = {Adv. Neural Inf. Process. Syst. 32},
	year      = {2019}
}

@inproceedings{huang2021streaming,
	title={Streaming $k$-{PCA}: Efficient guarantees for {Oja}'s algorithm, beyond rank-one updates},
	author={Huang, De and Niles-Weed, Jonathan and Ward, Rachel},
	booktitle={Proc. 34th Conf. Learn. Theory},
	pages={2463--2498},
	year={2021},
	organization={PMLR}
}

@inproceedings{jain2016streaming,
	title={Streaming {PCA}: Matching matrix {B}ernstein and near-optimal finite sample guarantees for {Oja}'s algorithm},
	author={Jain, Prateek and Jin, Chi and Kakade, Sham M and Netrapalli, Praneeth and Sidford, Aaron},
	booktitle={Proc. 29th Conf. Learn. Theory},
	pages={1147--1164},
	year={2016},
	organization={PMLR}
}

@article{oja1985stochastic,
	title={On stochastic approximation of the eigenvectors and eigenvalues of the expectation of a random matrix},
	author={Oja, Erkki and Karhunen, Juha},
	journal={J. Math. Anal. Appl.},
	volume={106},
	number={1},
	pages={69--84},
	year={1985},
	publisher={Elsevier}
}

@inproceedings{mitliagkas2013memory,
	title={Memory limited, streaming {PCA}},
	author={Mitliagkas, Ioannis and Caramanis, Constantine and Jain, Prateek},
	booktitle={Adv. Neural Inf. Process. Syst. 26},
	year={2013}
}

@inproceedings{shamir2016fast,
	title={Fast stochastic algorithms for {SVD} and {PCA}: Convergence properties and convexity},
	author={Shamir, Ohad},
	booktitle={Proc. 33rd Int. Conf. Mach. Learn.},
	pages={248--256},
	year={2016},
	organization={PMLR}
}

@inproceedings{hardt2014noisy,
	title={The noisy power method: A meta algorithm with applications},
	author={Hardt, Moritz and Price, Eric},
	booktitle={Adv. Neural Inf. Process. Syst. 27}, 
	year={2014}
}

@inproceedings{allen2017first,
	title={First efficient convergence for streaming $k$-{PCA}: A global, gap-free, and near-optimal rate},
	author={Allen-Zhu, Zeyuan and Li, Yuanzhi},
	booktitle={Proc. 2017 IEEE 58th Annu. Symp. Found. Comput. Sci.},
	pages={487--492},
	year={2017},
	organization={IEEE}
}

@article{li2018near,
	title={Near-optimal stochastic approximation for online principal component estimation},
	author={Li, Chris Junchi and Wang, Mengdi and Liu, Han and Zhang, Tong},
	journal={Math. Program.},
	volume={167},
	number={1},
	pages={75--97},
	year={2018},
	publisher={Springer}
}

@inproceedings{shamir2015stochastic,
	title={A stochastic {PCA} and {SVD} algorithm with an exponential convergence rate},
	author={Shamir, Ohad},
	booktitle={Proc. 32nd Int. Conf. Mach. Learn.},
	pages={144--152},
	year={2015},
	organization={PMLR}
}

@article{liang2023optimality,
	title={On the optimality of {Oja}'s algorithm for online {PCA}},
	author={Liang, Xin},
	journal={Stat. Comput.},
	volume={33},
	number={3},
	pages={62},
	year={2023},
	publisher={Springer}
}

@inproceedings{balsubramani2013streaming,
	title     = {Fast Convergence of Incremental {PCA}},
	author    = {Balsubramani, Akshay and Dasgupta, Sanjoy and Freund, Yoav},
	booktitle = {Adv. Neural Inf. Process. Syst. 26},
	year      = {2013}
}

@article{krasulina1969method,
	title   = {The Method of Stochastic Approximation in the Determination of the Least Eigenvalue of a Symmetric Matrix},
	author  = {Krasulina, T. P.},
	journal = {USSR Comput. Math. Math. Phys.},
	volume  = {9},
	number  = {6},
	pages   = {189--195},
	year    = {1969}
}

@article{oja1982simplified,
	title   = {Simplified Neuron Model as a Principal Component Analyzer},
	author  = {Oja, Erkki},
	journal = {J. Math. Biol.},
	volume  = {15},
	number  = {3},
	pages   = {267--273},
	year    = {1982}
}

@article{lanczos1950iteration,
	title   = {An Iteration Method for the Solution of the Eigenvalue Problem of Linear Differential and Integral Operators},
	author  = {Lanczos, Cornelius},
	journal = {J. Res. Natl. Bur. Stand.},
	volume  = {45},
	pages   = {255--282},
	year    = {1950}
}

@inproceedings{vogels2019powersgd,
	title     = {{PowerSGD}: Practical Low-Rank Gradient Compression for Distributed Optimization},
	author    = {Vogels, Thijs and Karimireddy, Sai Praneeth and Jaggi, Martin},
	booktitle = {Adv. Neural Inf. Process. Syst. 32},
	year      = {2019}
}

@InProceedings{zhao2024galore,
	title     = {{GaLore}: Memory-Efficient {LLM} Training by Gradient Low-Rank Projection},
	author    = {Zhao, Jiawei and Zhang, Zhenyu and Chen, Beidi and Wang, Zhangyang and Anandkumar, Anima and Tian, Yuandong},
	booktitle = {Proc. 41st Int. Conf. Mach. Learn.},
	pages     = {61121--61143},
	year      = {2024},
	publisher = {PMLR}
}

@article{halko2011finding,
	title   = {Finding Structure with Randomness: Probabilistic Algorithms for Constructing Approximate Matrix Decompositions},
	author  = {Halko, Nathan and Martinsson, Per-Gunnar and Tropp, Joel A.},
	journal = {SIAM Rev.},
	volume  = {53},
	number  = {2},
	pages   = {217--288},
	year    = {2011}
}

@article{eckart1936approximation,
	title   = {The Approximation of One Matrix by Another of Lower Rank},
	author  = {Eckart, Carl and Young, Gale},
	journal = {Psychometrika},
	volume  = {1},
	number  = {3},
	pages   = {211--218},
	year    = {1936}
}

@article{golub1970singular,
	title   = {Singular Value Decomposition and Least Squares Solutions},
	author  = {Golub, Gene H. and Reinsch, Christian},
	journal = {Numer. Math.},
	volume  = {14},
	number  = {5},
	pages   = {403--420},
	year    = {1970}
}

@incollection{jolliffe2005principal,
	title={Principal component analysis},
	author={Jolliffe, Ian},
	booktitle={Encyclopedia of Statistics in Behavioral Science},
	editor={Everitt, Brian S. and Howell, David C.},
	year={2005},
	publisher={John Wiley \& Sons},
	address={Chichester, UK}
}

@article{chang2011libsvm,
	title={{LIBSVM}: A Library for Support Vector Machines},
	author={Chang, Chih-Chung and Lin, Chih-Jen},
	journal={ACM Trans. Intell. Syst. Technol.},
	volume={2},
	number={3},
	pages={27:1--27:27},
	year={2011}
}

@book{saad2011numerical,
	title     = {Numerical Methods for Large Eigenvalue Problems},
	author    = {Saad, Yousef},
	edition   = {2nd},
	year      = {2011},
	publisher = {Society for Industrial and Applied Mathematics},
	address   = {Philadelphia, PA}
}

@article{wang2024dual,
	title={Tackling Arbitrarily Heterogeneous Data in Asynchronous Stochastic Gradient Descent Without Worker Scheduling},
	author={Wang, Xiaolu and Sun, Yuchang and Wai, Hoi-To and Zhang, Jun},
	journal={INFORMS J. Comput.},
	year={2026},
	note={Articles in Advance}
}

@book{CheegerEbin1975Comparison,
	title     = {Comparison Theorems in Riemannian Geometry},
	author    = {Cheeger, Jeff and Ebin, David G.},
	publisher = {North-Holland},
	year      = {1975}
}

@book{absil2008optimization,
	title={Optimization Algorithms on Matrix Manifolds},
	author={Absil, P-A and Mahony, Robert and Sepulchre, Rodolphe},
	year={2008},
	publisher={Princeton University Press}
}

@article{absil2004riemannian,
	title={Riemannian geometry of {Grassmann} manifolds with a view on algorithmic computation},
	author={Absil, P-A and Mahony, Robert and Sepulchre, Rodolphe},
	journal={Acta Appl. Math.},
	volume={80},
	number={2},
	pages={199--220},
	year={2004},
	publisher={Springer}
}

@article{hu2020brief,
	title={A Brief Introduction to Manifold Optimization},
	author={Hu, Jiang and Liu, Xin and Wen, Zai-Wen and Yuan, Ya-Xiang},
	journal={J. Oper. Res. Soc. China},
	volume={8},
	number={2},
	pages={199--248},
	year={2020}
}

@article{bendokat2024grassmann,
	title={A {Grassmann} manifold handbook: Basic geometry and computational aspects},
	author={Bendokat, Thomas and Zimmermann, Ralf and Absil, P-A},
	journal={Adv. Comput. Math.},
	volume={50},
	number={1},
	pages={6},
	year={2024},
	publisher={Springer}
}

@book{boumal2023introduction,
	title={An Introduction to Optimization on Smooth Manifolds},
	author={Boumal, Nicolas},
	year={2023},
	publisher={Cambridge University Press}
}

@inproceedings{wang2023incremental,
	title={Incremental aggregated {Riemannian} gradient method for distributed {PCA}},
	author={Wang, Xiaolu and Jiao, Yuchen and Wai, Hoi-To and Gu, Yuantao},
	booktitle={Proc. 26th Int. Conf. Artif. Intell. Stat.},
	pages={7492--7510},
	year={2023},
	organization={PMLR}
}

@article{peng2019nonconvex,
	title={Nonconvex proximal incremental aggregated gradient method with linear convergence},
	author={Peng, Wei and Zhang, Hui and Zhang, Xiaoya},
	journal={J. Optim. Theory Appl.},
	volume={183},
	number={1},
	pages={230--245},
	year={2019},
	publisher={Springer}
}

@article{alimisis2024geodesic,
	title={Geodesic convexity of the symmetric eigenvalue problem and convergence of steepest descent},
	author={Alimisis, Foivos and Vandereycken, Bart},
	journal={J. Optim. Theory Appl.},
	volume={203},
	number={1},
	pages={920--959},
	year={2024},
	publisher={Springer}
}

@book{golub2013matrix,
	title={Matrix Computations},
	author={Golub, Gene H and Van Loan, Charles F},
	edition={4th},
	year={2013},
	publisher={Johns Hopkins University Press},
	address={Baltimore, MD}
}

@unpublished{aytekin2016analysis,
	title={Analysis and implementation of an asynchronous optimization algorithm for the parameter server},
	author={Aytekin, Arda and Feyzmahdavian, Hamid Reza and Johansson, Mikael},
	note={arXiv:1610.05507},
	year={2016}
}

@mastersthesis{krizhevsky2009learning,
	author      = {Krizhevsky, Alex},
	title       = {Learning multiple layers of features from tiny images},
	school      = {University of Toronto},
	year        = {2009}
}

@article{assran2020advances,
	title={Advances in asynchronous parallel and distributed optimization},
	author={Assran, Mahmoud and Aytekin, Arda and Feyzmahdavian, Hamid Reza and Johansson, Mikael and Rabbat, Michael G},
	journal={Proc. IEEE},
	volume={108},
	number={11},
	pages={2013--2031},
	year={2020},
	publisher={IEEE}
}

@inproceedings{wang2023linear,
	title={Linear Speedup of Incremental Aggregated Gradient Methods on Streaming Data},
	author={Wang, Xiaolu and Jin, Cheng and Wai, Hoi-To and Gu, Yuantao},
	booktitle={Proc. 2023 62nd IEEE Conf. Decis. Control},
	pages={4314--4319},
	year={2023},
	organization={IEEE}
}

@book{bertsekas2015parallel,
	title={Parallel and Distributed Computation: Numerical Methods},
	author={Bertsekas, Dimitri and Tsitsiklis, John},
	year={2015},
	publisher={Athena Scientific}
}

@article{vanli2018global,
	title={Global convergence rate of proximal incremental aggregated gradient methods},
	author={Vanli, N Denizcan and Gurbuzbalaban, Mert and Ozdaglar, Asuman},
	journal={SIAM J. Optim.},
	volume={28},
	number={2},
	pages={1282--1300},
	year={2018},
	publisher={SIAM}
}

@article{zhang2021proximal,
	title={Proximal-like incremental aggregated gradient method with linear convergence under {Bregman} distance growth conditions},
	author={Zhang, Hui and Dai, Yu-Hong and Guo, Lei and Peng, Wei},
	journal={Math. Oper. Res.},
	volume={46},
	number={1},
	pages={61--81},
	year={2021},
	publisher={INFORMS}
}

@article{blatt2007convergent,
	title={A convergent incremental gradient method with a constant step size},
	author={Blatt, Doron and Hero, Alfred O and Gauchman, Hillel},
	journal={SIAM J. Optim.},
	volume={18},
	number={1},
	pages={29--51},
	year={2007},
	publisher={SIAM}
}

@article{gurbuzbalaban2017convergence,
	title={On the convergence rate of incremental aggregated gradient algorithms},
	author={Gurbuzbalaban, Mert and Ozdaglar, Asuman and Parrilo, Pablo A},
	journal={SIAM J. Optim.},
	volume={27},
	number={2},
	pages={1035--1048},
	year={2017},
	publisher={SIAM}
}

@article{higham1986computing,
	title={Computing the Polar Decomposition---with Applications},
	author={Higham, Nicholas J.},
	journal={SIAM J. Sci. Stat. Comput.},
	volume={7},
	number={4},
	pages={1160--1174},
	year={1986}
}

\end{document}